\documentclass[10pt]{amsart}

\usepackage{amssymb,eucal,mathabx}
\usepackage{amscd}
\usepackage[all,cmtip]{xy}
\usepackage{rotating}
\usepackage{mathrsfs,soul,stmaryrd,calrsfs}
\usepackage{upgreek}
\usepackage{tcolorbox}
\usepackage[bookmarksdepth=2,linktoc=page,colorlinks,linkcolor={red!70!black},citecolor={red!80!black},urlcolor={blue!80!black}]{hyperref}\newcommand{\arxiv}[1]{\href{http://arxiv.org/pdf/#1}{arXiv:#1}}
\usepackage[colorinlistoftodos,bordercolor=orange,backgroundcolor=orange!20,linecolor=orange,textsize=footnotesize]{todonotes} \makeatletter \providecommand \@dotsep{5} \def\listtodoname{List of Todos} \def\listoftodos{\@starttoc{tdo}\listtodoname} \makeatother 

\usepackage{amsmath,amssymb,epigraph,enumitem}
\usepackage{tikz}
\usepackage{tikz-cd}

\usepackage[top=2.5cm,bottom=2.5cm,right=3.35cm,left=3.35cm]{geometry}

\usepackage{etoolbox}
\makeatletter
\patchcmd{\@startsection}{\@afterindenttrue}{\@afterindentfalse}{}{}             
\makeatother

\newtheorem{theorem}{Theorem }[section]

\newtheorem{lemma}[theorem]{Lemma}

\newtheorem{remark}[theorem]{Remark}

\newtheorem{proposition}[theorem]{Proposition}

\theoremstyle{definition}
\newtheorem{df}[theorem]{Definition}
\newtheorem{rem}[theorem]{Remark}
\newtheorem{ex}[theorem]{Example}

\def\I{\mathbf{I}}

\def\eop{\hspace*{\fill}{\footnotesize$\blacksquare$}}

\setul{}{1.1pt}
\setulcolor{orange}

\newcommand{\mG}{\mathcal{G}}
\newcommand{\mC}{\mathcal{C}}
\newcommand{\mD}{\mathcal{D}}
\newcommand{\mQ}{\mathcal{Q}}

\newcommand{\wB}{\widehat{B}}
\newcommand{\wC}{\widehat{C}}
\usepackage{graphicx}
\usepackage{tikz}
\usepackage{graphicx}
\usepackage[all,cmtip]{xy}
\usepackage{rotating,multicol}
\usetikzlibrary{decorations.markings, shapes.misc, shadows}
\usetikzlibrary{matrix}
\usetikzlibrary{arrows}
\usepackage{pgf}
\usepackage{lscape}
\usepackage{listings}
\usetikzlibrary{positioning}
\usepackage{mdframed}
\newcommand*\circled[1]{\tikz[baseline=(char.base)]{
            \node[shape=circle,draw,inner sep=2pt] (char) {#1};}}

\newcommand{\mO}{\mathcal{O}}

\newcommand{\cB}{\mathcal{B}}

\newcommand{\cF}{\mathcal{F}}

\newcommand{\cN}{\mathcal{N}}
\newcommand{\cO}{\mathcal{O}}
\newcommand{\mP}{\mathcal{P}}

\newcommand{\mS}{\mathcal{S}}
\newcommand{\cT}{\mathcal{T}}

\newcommand{\cX}{\mathcal{X}}

\newcommand{\wt}{\widetilde}
\newcommand{\ol}{\overline}
\newcommand{\hT}{\mathbf{T}}

\newcommand{\bP}{\mathbb{P}}
\newcommand{\mL}{\mathcal{L}}

\DeclareSymbolFont{sfoperators}{OT1}{bch}{m}{n} \DeclareSymbolFontAlphabet{\mathsf}{sfoperators} \makeatletter\def\operator@font{\mathgroup\symsfoperators}\makeatother 
\DeclareMathOperator{\Bands}{Bands}
\DeclareMathOperator{\Crowds}{Crowds}
\DeclareMathOperator{\CrowdActivities}{CrowdActivities}
\DeclareMathOperator{\Groups}{Groups}
\DeclareMathOperator{\GroupActions}{GroupActions}
\DeclareMathOperator{\Rings}{Rings}
\DeclareMathOperator{\Sets}{Sets}
\DeclareMathOperator{\Alg}{Alg}
\DeclareMathOperator{\Sch}{Sch}
\DeclareMathOperator{\Spec}{\mathsf{Spec}}
\DeclareMathOperator{\Hom}{Hom}
\DeclareMathOperator{\Aut}{Aut}

\DeclareMathOperator{\SL}{SL}
\DeclareMathOperator{\PGL}{\mathsf{PGL}}
\DeclareMathOperator{\Gr}{\mathsf{Gr}}
\DeclareMathOperator{\Fl}{Fl}
\DeclareMathOperator{\Dr}{Dr}
\DeclareMathOperator{\sign}{sign}
\DeclareMathOperator{\rk}{rk}
\DeclareMathOperator{\ComRing}{ComRing}

\newcommand\A{{\mathbb A}}

\newcommand\F{{\mathbb F}}

\newcommand\bI{{\mathbf I}}

\newcommand\K{{\mathbb K}}

\newcommand\N{{\mathbb N}}

\renewcommand\P{{\mathbb P}}
\newcommand\Q{{\mathbb Q}}
\newcommand\R{{\mathbb R}}
\renewcommand\S{{\mathbb S}}
\newcommand\T{{\mathbb T}}

\newcommand\Z{{\mathbb Z}}

\newcommand\tA{{\mathsf A}}
\newcommand\tB{{\mathsf B}}
\newcommand\tC{{\mathsf C}}
\newcommand\tD{{\mathsf D}}
\newcommand\cG{{\mathcal G}}

\newcommand\cP{{\mathcal P}}
\newcommand\cR{{\mathcal R}}

\newcommand\cW{{\mathcal W}}
\renewcommand\1{{\mathbf 1}}
\newcommand\br{{\mathbf r}}
\renewcommand{\alpha}{\upalpha}
\renewcommand{\delta}{\updelta}
\renewcommand{\theta}{\uptheta}
\renewcommand{\tau}{\uptau}
\renewcommand{\rho}{\uprho}
\renewcommand{\varphi}{\upvarphi}
\renewcommand{\Gamma}{\Upgamma}
\renewcommand{\Delta}{\Updelta}
\renewcommand{\Pi}{\Uppi}
\renewcommand{\Omega}{\Upomega}
\renewcommand{\iota}{\upiota}
\renewcommand{\mu}{\upmu}

\renewcommand\emptyset\varnothing
\newcommand\Fun{{\F_1}}
\newcommand\Funpm{{\F_1^\pm}}
\newcommand\gen[1]{{\langle#1\rangle}}

\newcommand{\bgenquot}[2]{#1\hspace{-2.0pt}\sslash\hspace{-2.0pt}\gen{#2}}
\newcommand{\matrixtwo}[4]{\big[\begin{smallmatrix} #1 & #2 \\ #3 & #4\end{smallmatrix}\big]}
\newcommand{\matrixthree}[9]{\Big[\begin{smallmatrix} #1 & #2 & #3 \\ #4 & #5 & #6 \\ #7 & #8 & #9 \end{smallmatrix}\Big]}
\newcommand{\tinymatrix}[6]{\Big[\begin{smallmatrix} #1 & #4 \\ #2 & #5 \\ #3 & #6 \end{smallmatrix}\Big]}
\newcommand{\tinyvector}[3]{{}\Big[\begin{smallmatrix} #1 \\ #2 \\ #3 \end{smallmatrix}\Big]}

\title[On band schemes, crowds and $\F_1$-structures]{Towards the horizons of Tits's vision | on band schemes, crowds and $\F_1$-structures}

\subjclass[2000]{}

\author{Oliver Lorscheid}
\address{Oliver Lorscheid, University of Groningen, Nijenborgh 9, 9747 AG Groningen, the Netherlands, and IMPA, Estrada Dona Castorina, CEP 22460-320, Rio de Janeiro, Brazil}
\email{oliver@impa.br}

\author{Koen Thas}

\address{Koen Thas, {Ghent University},
{Department of Mathematics: Algebra and Geometry},
{Krijgslaan 281, S25, B-9000 Ghent, Belgium}}
\email{koen.thas@gmail.com}
\thanks{}
\date{}

\begin{document}

\newcommand{\epiarrow}[1]{%
\parbox{#1}{\tikz{\draw[thick,->>](0,0)--(#1,0);}}
}

\newcommand{\monoarrow}[1]{%
\parbox{#1}{\tikz{\draw[thick,right hook->](0,0)--(#1,0);}}
}

\begin{abstract}
 This text is dedicated to Jacques Tits's ideas on geometry over $\Fun$, the field with one element. In a first part, we explain how thin Tits geometries surface as rational point sets over the Krasner hyperfield, which links these ideas to combinatorial flag varieties in the sense of Borovik, Gelfand and White and $\Fun$-geometry in the sense of Connes and Consani. A completely novel feature is our approach to algebraic groups over $\Fun$ in terms of an alteration of the very concept of a group. In the second part, we study an incidence-geometrical counterpart of (epimoprhisms to) thin Tits geometries; we introduce and classify all $\F_1$-structures on $3$-dimensional projective spaces over finite fields. This extends recent work of Thas and Thas \cite{part4} on epimorphisms of projective planes (and other rank $2$ buildings) to thin planes.   
\end{abstract}
 
\maketitle

\begin{quote}\footnotesize 
\textsc{Epigraph.} 13. \textit{Les groupes de Chevalley sur les ``corps de caract\'eristique 1.''} Nous avons vu au n${}^{\text{o}}$ 9 que les groupes de Chevalley sur un corps donn\'e $K$ et les g\'eom\'etries sur $K$ correspondant \`a tous les sch\'emas de Witt-Dynkin sont d\'etermin\'es, par l'interm\'ediaire des propositions g\'en\'erales des n${}^{\text{os}}$ 5 \`a 8, d\`es qu'on connait les g\'eom\'etries correspondant aux sch\'emas de la fig.\ 4. On peut alors songer \`a associer \`a ces derniers d'autres g\'eom\'etries que celles indiqu\'ees au n${}^{\text{o}}$ 9, et \`a rechercher si les propositions g\'en\'erales des n${}^{\text{os}}$ 5 \`a 8 conduisent encore \`a associer aux autres sch\'emas de Witt-Dynkin (ou \'eventuellement, \`a certains d'entre eux) des g\'eom\'etries univoquement d\'etermin\'ees. C'est ce que nos ferons ici.

Nous designerons par $K=K_1$ le ``corps de carat\'eristique $1$'' form\'e du seul \'el\'ement $1=0$. Il est naturel d'appeler \textit{espace projectif \`a n dimensions sur $K$}, un ensemble $\cP_{n+1}$ de $n+1$ points dont tous les sous-ensembles sont consid\'er\'es comme des vari\'et\'es lin\'eaires, la dimension d'une vari\'et\'e \'etant le nombre de points qui la constituent diminu\'e d'une unit\'e, et \textit{projectivit\'e de $\cP_n$}, une permutation quelconque de ces points. On d\'efinit alors, suivant (3.1), la \textit{g\'eometrie projective \`a n dimensions sur $K$}, $\Pi_{n,K}$. 

\hfill \textit{Jacques Tits, excerpt from \cite{Tits57}}%
\end{quote}

\begin{tcolorbox}
 \begin{small}
  \setcounter{tocdepth}{1}
  \tableofcontents
 \end{small}
\end{tcolorbox}

\section{Introduction}

The idea of a field with one element was first perceived by Jacques Tits in the 1950s who saw a parallel between incidence geometries that stem from algebraic geometry over finite fields and those that stem from combinatorial group theory. He communicated his ideas in several talks and included a short section in his paper \cite{Tits57} (as reproduced in the epigraph), which is the only published account of what he had in mind.

The central theme of Tits's paper is a framework for incidence geometries that captures both worlds | geometry over Galois fields $\F_q$ with $q$ elements and combinatorial group theory | which behaves analogously, but with ``parameter'' $q=1$. In particular, Tits asks for an algebro-geometric explanation of these combinatorial geometries in terms of algebraic geometry over the elusive field $\Fun$ with one element.

\subsection*{Tits geometries}

We begin with the revision of some key ideas from \cite{Tits57}. A \emph{Tits geometry}\footnote{As the reader might tell, Tits did not use a terminology that includes his own name. He rather used the generic word ``geometry'' for what we call ``Tits geometry'' in his honour.} $\Gamma=\Gamma(E;F_1,\dotsc,F_n;\bI;G)$ consists of a set $E$ together with a partition
\[
 E \ = \ \bigsqcup_{i=1}^n \ F_i
\]
into disjoint subsets, an incidence relation $\bI \subset E\times E$ and its automorphism group $G=\Aut(\Gamma)$, which consists of all permutations of $E$ that preserve the partition and the incidence relation. 

The leading example of a Tits geometry is that of the collection $E$ of all non-trivial and proper $K$-linear subspaces of $K^{n}$ for a field $K$ with $F_d$ being the family of $d$-dimensional subspaces and $\bI$ being the symmetrisized inclusion relation of subspaces. 
This Tits geometry is linked intimately to projective algebraic geometry over $K$: the elements of $F_d$ correspond to the $(d-1)$-dimensional linear subspaces of the projective space $\P^{n-1}(K)$ or, equivalently, to the points of the Grassmannian $\Gr(d,n)(K)$. This defines the \emph{projective Tits geometry $\Pi_{n,K} = \mathbb{P}^{n - 1}(K)$ of dimension $n-1$ over $K$}.

Similarly all other (simple) Chevalley groups (of classical types $\tA$, $\tB$, $\tC$ and $\tD$) appear as (simple) automorphism groups of Tits geometries. The notion of Tits geometry has later been made obsolete by the invention of buildings (\cite{Bruhat-Tits72}, \cite{Tits95}), which covers the geometry of Chevalley groups and other combinatorial geometries, such as trees and Coxeter complexes, in an effective and powerful way.

\subsection*{Geometry over \texorpdfstring{$\Fun$}{F1}}

In the example of the projective Tits geometry $\Pi_{n,K}$, we find a uniform behaviour for finite fields $K=\F_q$. For example, an element $W\in F_d$ is incident to exactly $[d]_q=1+q+\dotsb+q^{d-1}$ elements in $F_{d-1}$ and to $[n-d]_q$ elements in $F_{d+1}$.

The Coxeter complex of $S_n$, considered as a Tits geometry, behaves like the limit $q\to1$ of $\Pi_{n,\F_q}$: its set $E$ consists of all non-empty and proper subsets of $\{0,\dotsc,n\}$, which is partitioned into the collections $F_d$ of subsets with $d-1$ elements, together with the inclusion relation $i$. Its automorphism group is the permutation group $S_n$. In this case, an element $W\in F_d$ is incident to $[d]_1=1^0+\dotsb+1^{d-1}=d$ elements in $F_{d-1}$. For this reason, we denote the Coxeter complex of $S_n$ by $\Pi_{n,\Fun}$.

Rephrasing Tits's words, we would like to find an explanation of this Tits geometry in terms of an algebraic geometry over $\Fun$ that features models of Grassmannians $\Gr(d,n)$ and the projective linear group $\PGL(n)$ such that $F_d=\Gr(d,n)(\Fun)$ and $S_n=\PGL(n,\Fun)$.

\subsection*{Influence on later developments}

Tits's arcane remark from the epigraph has woken a serious interest by other researchers only around 30 years later when the concept of a field with one element was interwoven with other ideas, such as approaches to the Riemann Hypothesis (\cite{Manin95}) and the ABC Conjecture (\cite{Smirnov92}); cf.\ \cite{LopezPena-Lorscheid11,Lorscheid18} for more detailed expositions.

Since then, it has attracted much attention (\cite{Borger09,Cohn04,Connes-Consani10a,Connes-Consani11a,Deitmar05,Deitmar13,Kapranov-Smirnov,Kurokawa05,LopezPena-Lorscheid11,Manin95,Soule04,Toen-Vaquie09}) and found several solutions (\cite{Lorscheid12b,Lorscheid16,Lorscheid18b}). In this text, we present yet another solution that ties in with the $\Fun$-perspective on hyperrings (cf.\ \cite{Connes-Consani10b}) and which might shed new light on algebraic groups in matroid theory and tropical geometry.

\subsection*{Combinatorial flag varieties after Borovik, Gelfand and White}

Tits's dream of algebraic groups over the ``field of characteristic one'' surfaces as the outer layer of the deep and mysterious geometry of the combinatorial flag varieties $\Omega_W$ of Gelfand-Stone-Rybnikov, which are simplicial complexes whose simplices correspond to flag matroids. This insight\footnote{We like to remark at this point that such Coxeter complexes have in general a rich structure, opposed to the cited passage by Gelfand-Stone-Rybnikov. For example, the theory of finite generalized octagons of order $(1,s)$ is essentially equivalent to the theory of finite generealized quadrangles of order $(s,s)$. (We also refer to part 2 of this paper for a more detailed discussion.)} can be found in \cite[section 7.14]{Borovik-Gelfand-White03}:
\begin{quote}\small
 Many geometries over fields have formal analogues which can be thought of as geometries over the field of $1$ element. For example, the projective plane over the field $\F_q$ has $q^2+q+1$ points and the same number of lines; every line in the plane has $q+1$ points. When $q=1$, we have a plane with three points and three lines, i.e., a triangle. The flag complex of the triangle is a thin building of type $\tA_2 = \mathsf{Sym}_3$. In general, the Coxeter complex $\cW$ of a Coxeter group $W$ is a thin building of type $W$ and behaves like the building of type $W$ over the field of $1$ element.

 However, the Coxeter complex has a relatively poor structure. In many aspects, $\Omega_W$ and $\Omega_W^\ast$ are more suitable candidates for the role of a ``universal'' combinatorial geometry of type $W$ over the field of $1$ element.
\end{quote}
The combinatorial flag variety $\Omega_{S_n}$ is the simplicial complex whose vertices are matroids on $E=\{1,\dotsc,n\}$ and whose higher simplices correspond to flags of matroids on $E$. The Coxeter complex of $S_n$ embeds into the combinatorial flag variety $\Omega_{S_n}$ as the subcomplex whose vertices are matroids on $\{1,\dotsc,n\}$ that decompose into a direct sum of matroids on $1$ element.

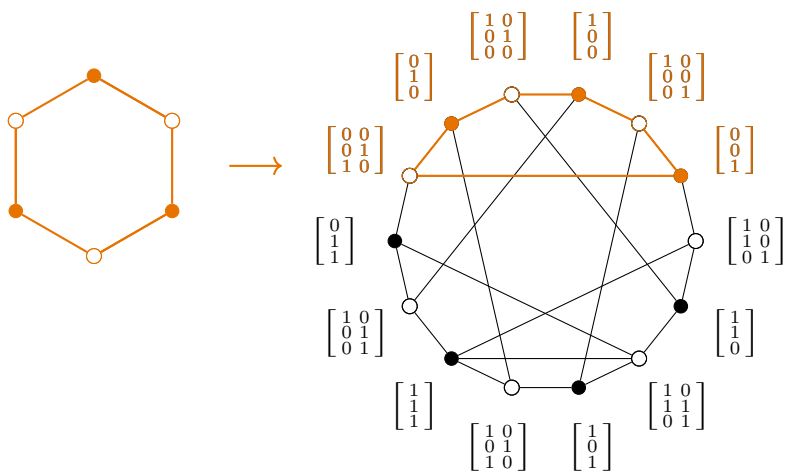
\begin{figure}[htb]
 \[
  \newcounter{tikz-counter}
  \begin{tikzpicture}
   \node (origin) {};
   \foreach \a in {1,...,14}{\draw (\a*360/7+0/7: 2cm) node [draw,circle,inner sep=2pt] (l\a) {};}
   \foreach \a in {1,...,7}{\draw (\a*360/7+180/7: 2cm) node [draw,circle,inner sep=1.8pt,fill] (p\a) {};
                              \draw [-] (p\a) -- (l\a);
                              \setcounter{tikz-counter}{\a};
                              \addtocounter{tikz-counter}{1};
                              \draw [-] (p\a) -- (l\arabic{tikz-counter});
                              \addtocounter{tikz-counter}{2};
                              \draw [-] (p\a) -- (l\arabic{tikz-counter});
                             }
   \draw [-] (p4) -- (l6);                          
   \node at      (0:2.8cm) {$\tinymatrix 110001$};                          
   \node at  (360/7:2.8cm) {$\tinymatrix 100001$};                          
   \node at  (720/7:2.8cm) {$\tinymatrix 100010$};                          
   \node at (1080/7:2.8cm) {$\tinymatrix 001010$};                          
   \node at (1440/7:2.8cm) {$\tinymatrix 100011$};                          
   \node at (1800/7:2.8cm) {$\tinymatrix 101010$};                          
   \node at (2160/7:2.8cm) {$\tinymatrix 110011$};                          
   \node at  (180/7:2.8cm) {$\tinyvector 001$};                          
   \node at  (540/7:2.8cm) {$\tinyvector 100$};                          
   \node at  (900/7:2.8cm) {$\tinyvector 010$};                          
   \node at (1260/7:2.8cm) {$\tinyvector 011$};                          
   \node at (1620/7:2.8cm) {$\tinyvector 111$};                          
   \node at (1980/7:2.8cm) {$\tinyvector 101$};                          
   \node at (2340/7:2.8cm) {$\tinyvector 110$};   
   \node[color=orange!90!black] at  (360/7:2.8cm) {$\tinymatrix 100001$};                          
   \node[color=orange!90!black] at  (720/7:2.8cm) {$\tinymatrix 100010$};                          
   \node[color=orange!90!black] at (1080/7:2.8cm) {$\tinymatrix 001010$};                          
   \node[color=orange!90!black] at  (180/7:2.8cm) {$\tinyvector 001$};                          
   \node[color=orange!90!black] at  (540/7:2.8cm) {$\tinyvector 100$};                          
   \node[color=orange!90!black] at  (900/7:2.8cm) {$\tinyvector 010$};                          
   \foreach \a in {1,...,3}{\draw (\a*360/7+0/7: 2cm) node [draw,color=orange!90!black,circle,inner sep=2pt] (l\a) {};}
   \foreach \a in {0,...,2}{\draw (\a*360/7+180/7: 2cm) node [draw,color=orange!90!black,circle,inner sep=1.8pt,fill=orange!90!black] (p\a) {};}
   \draw [-,color=orange!90!black,thick] (p0) -- (l1) -- (p1) -- (l2) -- (p2) -- (l3) -- (p0);
   \foreach \a in {1,...,6}{\draw [xshift=-6cm,yshift=1cm] (\a*360/3+30: 1.2cm) node [draw,color=orange!90!black,circle,inner sep=2pt] (l\a) {};}
   \foreach \a in {1,...,3}{\draw [xshift=-6cm,yshift=1cm] (\a*360/3+90: 1.2cm) node [draw,color=orange!90!black,circle,inner sep=1.8pt,fill=orange!90!black] (p\a) {};
                              \draw [-,color=orange!90!black,thick] (p\a) -- (l\a);
                              \setcounter{tikz-counter}{\a};
                              \addtocounter{tikz-counter}{1};
                              \draw [-,color=orange!90!black,thick] (p\a) -- (l\arabic{tikz-counter});
                              \addtocounter{tikz-counter}{2};
                              \draw [-,color=orange!90!black,thick] (p\a) -- (l\arabic{tikz-counter});
                             }
   \draw[->,color=orange!90!black,thick] (-4.2,1) -- (-3.5,1);                          
  \end{tikzpicture}
 \] 
 \caption{Embedding of the Coxeter complex of $S_3$ into $\Omega_{S_3}$.}
 \label{fig: Coxeter complex and combinatorial flag variety for S3}
\end{figure}

Figure \ref{fig: Coxeter complex and combinatorial flag variety for S3} illustrates the Coxeter complex of $S_3$ as a subspace of the combinatorial flag variety $\Omega_{S_3}$. The vertices are labelled by matrices that represent the corresponding matroids. 

\subsection*{Combinatorial flag varieties as rational point sets of \texorpdfstring{$\Fun$}{F1}-schemes}

The joint works of the second author with Baker (\cite{Baker-Lorscheid21}) and Jarra (\cite{Jarra-Lorscheid22}) on moduli spaces of (flag) matroids describe models of Grassmannians and flag varieties over the field with one element, or, more accurately, over the \emph{regular partial field $\Funpm$}, which can be thought of as a quadratic extension of $\Fun$ that adjoins an element $-1$.

These varieties are, in essence, closed subschemes of a product of projective spaces over $\Funpm$ that is defined by usual quadratic Pl\"ucker and incidence relations. We explain this theory in terms of \emph{band schemes}, which the second author introduces in joint work with Baker and Jin (\cite{Baker-Jin-Lorscheid}) as a simplification of ordered blue schemes that is suitable to cover many interesting aspects of $\Fun$-geometry. 

\subsection*{Band schemes}

To give a brief impression on band schemes: a {\em band} is a multiplicatively written commutative monoid $B$ with neutral element $1$ and absorbing element $0$, together with a \emph{null set} $N_B$, which is an ideal of the semiring $\N[B]/\gen{0}=\{\sum n_aa\mid n_a\in\N,\, a\in B-\{0\}\}$ (where only finite sums are considered) such that every $a\in B$ has a unique \emph{additive inverse} $-a\in B$ that satisfies $a+(-a)\in N_B$. 

For example, a ring $R$ is a band with null set $N_R=\{\sum n_aa\mid \sum n_aa=0\text{ in }R\}$. Two important examples of bands that are not rings are 
\begin{itemize}
 \item the \emph{regular partial field} $\Funpm=\{0,1,-1\}$ with null set $N_\Funpm=\{n_11+n_{-1}(-1)\mid n_1=n_{-1}\}$;
 \item the \emph{Krasner hyperfield} $\K=\{0,1\}$ with null set $N_\K = \{n_1\cdot 1\mid n_1\neq 1\}$.
\end{itemize}

Bands come with a suitable notion of morphism, which defines the category $\Bands$. A \emph{band scheme} is a functor $\cX:\Bands\to\Sets$ that is locally representable in a suitable sense. Examples of band schemes are representable, or \emph{affine}, band schemes $\Hom(C,-)$ where $C$ is a band. 

Another example is the \emph{projective $n$-space $\P^n$} with
\[
 \P^n(B) \ = \ \big\{ \, (a_0,\dotsc,a_n)\in B^{n+1} \,\big|\, a_i\in B^\times\text{ for some }i \, \big\} \; / \; B^\times
\]
where the unit group $B^\times=\{a\in B\mid ab=1\text{ for some }b\in B\}$ acts diagonally on $B^{n+1}$ (and $B$ is a band). 

Grassmannians $\Gr(r,n)$ and flag varieties $\Fl(\br;n)$ of type $\br=(r_1,\dotsc,r_s)$ are defined by the usual Pl\"ucker and incidence relations as subfunctors of (products of) projective spaces. They come with coordinate projections $\Fl(\br,n)\to\Fl(\br',n)$ whenever $\br'$ is a subtype of $\br$.

The band schemes $\Fl(\br,n)$ extend the usual flag varieties in the sense that $\Fl(\br,n)(k)$ agrees with the usual flag varieties over a field $k$. In particular, we recover the description of spherical buildings of type $\tA_{n-1}$ in terms of flag varieties over $\F_q$: the simplices of the building stay in bijection with 
 \[
  \coprod_{\br\in\Theta} \ \Fl(\br,n)(\F_q) \quad \text{where} \quad \Theta \ = \ \big\{ \, (r_1,\dotsc,r_s) \, \big| \, s>0,\; 0<r_1<\dotsb<r_s<n \, \big\},
 \]
and the face relations of simplices is given by the projections $\Fl(\br,n)(\F_q)\to\Fl(\br',n)(\F_q)$.

Taking $\K$-rational points recovers in an analogous fashion the combinatorial flag variety $\Omega_{S_n}$. 

\begin{theorem}\label{thmA}
 The simplices of the combinatorial flag variety $\Omega_{S_n}$ correspond bijectively to  
 \[
  \coprod_{\br\in\Theta} \ \Fl(\br,n)(\K) 
 \]
 and the face relation of simplices is given by the projections $\Fl(\br,n)(\K)\to\Fl(\br',n)(\K)$.
\end{theorem}

\subsection*{Crowds}

The aspect that is completely novel in this paper is the realization, and extension, of the group action of $S_n$ on its Coxeter complex and its combinatorial flag variety in terms of an action of band schemes. This theory starts with the generalization of groups to ``crowds.''

A {\em crowd} is a set $G$ with an \emph{identity} $1\in G$ and a \emph{crowd law} $R\subset G^3$, which is assumed to satisfy some suitable axioms that we explain in section \ref{subsection: definition of crowds}. A group $G$ is a crowd with crowd law $R=\{(a,b,c)\in G^3\mid abc=1\}$. Note that we can recover the group multiplication from the crowd law: the product $ab$ of two elements $a,b\in G$ is the unique element $c\in G$ for which there is a $\bar c\in G$ with $(a,b,\bar c),\, (c,\bar c,1)\in R$.

A natural notion of morphism turns crowds into a category $\Crowds$. This allows for an extension of algebraic groups to band schemes: we define an \emph{algebraic crowd} as a band scheme $\cG:\Bands\to\Crowds$ for which the crowd law $\cR:\Bands\to\Sets$ is also a band scheme.

The main example of this text is the algebraic crowd $\SL_n$ with
\[
 \SL_n(B) \ = \ \big\{ (a_{i,j})\in B^{n^2} \mid \det(a_{i,j})-1\in N_B \big\}
\]
where we write $a=(a_{i,j})_{i,j=1,\dotsc,n}$, whose crowd law $\cR(B)$ consists of all triples $(a,b,c)\in\SL_n(B)^3$ for which 
\[
 \sum_{k,l=1}^n \ a_{i,k} b_{k,l} c_{l,j} - \delta_{i,j} \ \in \ N_B
\]
for all $i,j=1,\dotsc,n$, as well as for all cyclic permutations of $(a,b,c)$. Note that the crowd $\SL_n(B)$ is in general not a group, which underlines the role of crowds for this theory.

\subsection*{Crowd activities}
Replacing the action $\theta:G\times X\to X$ of a group $G$ on a set $X$ by its graph $T=\{(a,x,y)\in G\times X\times X\mid a. x=y\}$ leads to the notion of a \emph{crowd activity} of a crowd $G$ on a set $X$ as a subset $T\subset G\times X\times X$. An \emph{algebraic crowd activity} of an algebraic crowd $\cG$ on a band scheme $\cX$ is a subscheme $\cT$ of $\cG\times\cX\times\cX$.

The classical action of the special linear group on flag varieties extends in a natural way to a crowd activity of the algebraic crowd $\SL_n$ on the band scheme $\Fl(\br,n)$. Taking $\K$-rational points exhibits a new type of symmetries for combinatorial flag varieties. More precisely, we have:

\begin{theorem}\label{thmB}
 The crowd $\SL_n(\K)$ contains the symmetric group $S_n$ as the subcrowd of all permutation matrices. The crowd activity of $\SL_n(\K)$ on $\Fl(\br,n)(\K)$ for various $\br$ defines a crowd activity on $\Omega_{S_n}$ whose restriction to $S_n$ is the usual action on $\Omega_{S_n}$.
\end{theorem}


We consider it as a problem of high interest to investigate the properties of this crowd action, which might be a suitable tool to tackle various combinatorial conjectures with a proven analogue in finite field geometry. Another interesting task is to extend the insights of this text to algebraic groups and flag varieties of other Dynkin types.

\subsection*{Epimorphisms and $\Fun$-structures}

Inspired by work done on covers of generalized quadrangles \cite{part1,part2}, the same authors started a systematic study of 
epimorphisms $\upepsilon$ of generalized polygons, both finite and infinite, to thin generalized polygons of the same gonality | in an obvious 
sense  $\Fun$-versions of the former. In particular, in \cite{part3}, Thas and Thas classified all epimorphisms 
\begin{equation}
\upepsilon: \ \Upgamma\ \epiarrow{1.1cm}\ \Updelta = K(3),    
\end{equation}
where $\Upgamma$ is a finite axiomatic (so not necessarily Desarguesian) projective plane, and $\Updelta$ is isomorphic to the complete graph on 
three vertices | that is, a projective plane of order $1$ (the $\Fun$-plane). In \cite{part3,part4}  it is also shown that such a classification 
cannot work for infinite ploygons, since free constructions of epimorphisms to thin polygons of the same gonality exist. In the second part of the present paper, we want to take 
the planar result a step | that is, one dimension | further, and we classify the epimorphisms 
\begin{equation}
\upepsilon: \ \mathbb{P}^3(\F_q)\ \epiarrow{1cm}\ \Updelta = K(4).
\end{equation}

Each such epimorphisms naturally endows $\mathbb{P}^3(\F_q)$ with an {\em $\Fun$-structure}. 

Note that 
we automatically work with projective spaces over division rings (and hence in the finite case, fields), since an old result of Veblen-Young tells 
us that axiomatic projective spaces of dimension at least three are defined over division rings anyhow. 

We summarize the main result of the second part as follows (and we refer to the appropriate sections for a defnition of type {\bf II}-planes):

\begin{theorem}
Let $k$ be a finite field.
Consider the surjective morphism $\upepsilon: \ \mP = \mathbb{P}^3(k) \ \epiarrow{1cm} \ K(4)$; then the possible $\Fun$-structures endowed by $\upepsilon$ on $\mathbb{P}^3(k)$ are up to 
a permutation of $A, B, C, D$, described as follows.
 \begin{quote}
\begin{mdframed}[linewidth=2.5,linecolor=orange, topline=false,rightline=false,bottomline=false]
{\em There is a line $U$ such that $\vert U \cap A \vert \geq 1$ and $\vert U \cap B \vert \geq 1$, and $U \subseteq A \cup B$, and we distinguish two  cases.}
\item[{\rm $\circled{i}$}]
There is one plane $\Pi_C$ of Type {\bf II} containing $U$, and its unique point not in $A \cup B$ is contained in $C$; for any other plane $\Pi$ on $U$ we have that the points 
of $\Pi \setminus U$ are contained in $D$;
\item[{\rm $\circled{ii}$}] 
There is a line $\widehat{U}$ for which $\vert C \cap \widehat{U} \vert \geq 1$, $\vert D \cap \widehat{U} \vert \geq 1$ and $C \cup D = \widehat{U}$; each point $u$ which is not 
incident with $U$ nor $\widehat{U}$ is incident with a unique line $W$ which meets $U$ and $\widehat{U}$, and if $V \cap U$ is a point of $A$, respectively $B$, then 
all points of $V \setminus V \cap \widehat{U}$ are points of $A$, respectively $B$.
\end{mdframed} 
\end{quote}
\eop 
\end{theorem}

\bigskip
\subsection*{In short}

In this paper, we study the diagram below. In the diagram, we consider geometric objects $\mG = \widehat{\Spec}(A)$ which are defined 
over $\Fun$ through some (to-be-defined) scheme theory $\widehat{\Spec}$ (where we only mention the affine part for the sake of convenience), where 
$A$ is an object in a category which generalizes the category of commutative rings. 

Part 1 considers the first line of the diagram (the algebro-geometric part) on the level of algebraic groups, and focuses on the special linear groups $\SL(n,\cdot)$, while in part 2 we study its incidence-geometrical counterpart, which 
is presented by the second line, focusing on the natural geometric modules on which $\SL(n,\cdot)$ acts, namely projective spaces. (The second line 
represents the geometry of the rational points.) \\

\bigskip
\begin{center}
{\scalebox{1.1}{
\begin{tikzpicture}[scale=2.2,text height=1.5ex, text depth=0.25ex]
\node (a0) at (0,3) {$\mG = \widehat{\Spec}(A)$};
\node (a) [right=of a0] {};
\node (a1) [right=of a] {$\Updelta = \widehat{\Spec}(\Fun)$};
\node (a2) [right=of a1,scale=1.2,draw,color=orange,fill=white,
  shape=rounded rectangle,thick,
  drop shadow={opacity=.5,shadow xshift=0pt},
  minimum width=1cm]
  {part 1};

\node (b0) [below=of a0] {$\mG(k)$};
\node (b1) [below=of a1] {$\Updelta(\Fun)$};
\node (b2) [right=of b1,scale=1.2,draw,color=orange,fill=white,
  shape=rounded rectangle,thick,
  drop shadow={opacity=.5,shadow xshift=0pt},
  minimum width=1cm]
  {part 2};

\draw[->>,font=\scriptsize,thick]
(a0) edge node[auto] {} (a1)
(b0) edge node[auto] {} (b1);

\draw[->,font=\scriptsize,thick]
(a0) edge node[auto] {} (b0)
(a1) edge node[auto] {} (b1);

\draw[->,font=\scriptsize,orange,thick]
(a1) edge node[auto] {} (a2)
(b1) edge node[auto] {} (b2);

\end{tikzpicture}
}}
\end{center}

\vspace*{0.5cm}
\ 
\part*{Part $0$. A very short casual introduction to the definition of $\F_1$} 
\addcontentsline{toc}{subsection}{Part $0$. A very short casual introduction to the definition of $\F_1$}

\section{The field \texorpdfstring{$\Fun$}{F1}}

In this section we state some algebraic definitions and properties of $\Fun$. 

\subsection*{The ``field'' \texorpdfstring{$\Fun$}{F1}}

One way of defining $\Fun$ | to fix ideas | is as the set 
\begin{equation}
\Big(\{ 0, 1\}, \cdot\Big),
\end{equation}
where only the classical multiplication is present (so one does not allow addition). 
Although this definition is very simple, we will show in the next subsection that by relaxing (or better: releasing) the addition of $\F_2$ to obtain $\Fun$, its absolute galois group already is an object of interest. 

First note that the algebraic closure $\overline{\Fun}$ of $\Fun$ is the group of all roots of unity (which is isomorphic to $(\Q/\Z,+)$); this is because the only polynomial equations we consider in this simple setting are of the form $x^m = 1$ (with $m$ a positive integer). Another way to look at it goes as follows: let $p$ be any prime, and consider the algebraic closure $\overline{\F_p}$ of $\F_p$; it consists of all finite field extensions $\F_{p^n}$ of $\F_p$. If we represent $\Fun$ as above, one standard way to define $\F_{1^n}$, with $n$ a positive integer, is as the set   
\begin{equation}
\Big(\{ 0\} \cup \upmu_n, \cdot\Big),
\end{equation}
where $\upmu_n$ is the cyclic group of order $n$. Now let $n$ vary to obtain $\overline{\Fun}$. And note that this approach is independent of the chosen prime $p$.

\subsection*{Absolute Galois group of \texorpdfstring{$\Fun$}{F1}}

We can now easily calculate the absolute Galois group 
\begin{equation}
\mathsf{Gal}(\overline{\Fun}/\Fun) \cong \Aut(\overline{\Fun}) 
\end{equation}
of $\Fun$. First note that $(\Q/\Z,+)$ is the direct sum of its so-called {\em Pr\"{u}fer $p$-groups}\index{Pr\"{u}fer $p$-group} $\Z(p^\infty)$, where $p$ ranges over the primes. (One can define $\Z(p^\infty)$ simply as the Sylow $p$-subgroup of $\Q/\Z$; alternatively, it is the group of $p^n$-roots of unity, where $n \geq 1$ runs over the positive integers.) An endomorphism of $\Q/\Z$ induces an endomorphism in each component $\Z(p^\infty)$ and vice versa, each choice of a set of endomorphisms ${(\upalpha_p)}_p$ (where $\upalpha_p$ is an endomorphism of $\Z(p^\infty)$) as $p$ ranges over the primes, defines an endomorphism of $\Q/\Z$. So the ring of endomorphisms of $\Q/\Z$ is isomorphic to  the direct product of the endomorphism rings $\mathsf{End}(\Z(p^\infty))$ of the $\Z(p^\infty)$, as $p$ ranges over the primes. 
It is well known that $\mathsf{End}(\Z(p^\infty))$ is isomorphic to the ring of $p$-adic integers $\Z_p$. Passing to units, we finally obtain that the automorphism group we are seeking is given by $\prod_p\Z_p^\times$, which in its turn is isomorphic to the group $\widehat{\Z}^\times$, the group of units in the ring of profinite integers $\widehat{\Z}$. 

Note that $(\widehat{\Z},+)$ is isomorphic to $\mathsf{Gal}(\overline{\F_q}/\F_q) \cong \Aut(\overline{\F}_q)$, where $q$ is any prime power. It is also important to note that since 
\begin{equation}
(\widehat{\Z},+) \ \cong\ \mathsf{Gal}(\overline{\F_q}/\F_q)
\end{equation}  
is independent of the choice of $q$, the profinite integers capture a lot of the $\Fun$-behavior underlying finite fields.

\ 
\part{Algebraic groups over \texorpdfstring{$\Fun$}{F1}}

The vessel for our realization to Tits's vision on algebraic groups over $\Fun$ is the concept of band schemes. For the versed reader, band schemes are particular types of ordered blue schemes. But we will give an independent non-technical account that allows the new-comer to $\Fun$-geometry follow this text. 

Our approach is example driven. Conceptually this theory promises to lead to a new interesting perspective on algebraic groups over $\Fun$ and over other exotic structures, such as the tropical numbers. In so far, we stimulate future research on algebraic groups over $\Fun$ in the proposed language.

\section{Bands}
\label{section: bands}

Bands and band schemes are introduced in the forthcoming paper \cite{Baker-Jin-Lorscheid} by Baker, Jin and the first author. They provide a simplified account to $\Fun$-geometry (opposed to ordered blue schemes) that is tailored to certain central applications of $\Fun$-geometry. In particular, band schemes serve for the purposes of this text.

\subsection{Bands}
\label{subsection: bands}

A \emph{pointed monoid} is a (multiplicatively written) commutative monoid $A$ with neutral element $1$ together with an absorbing element $0$ that satisfies $0\cdot a=0$ for all $a\in A$. We can embed $A$ as a multiplicative submonoid in the semiring
\[\textstyle
 A^+ \ = \ \N[A]/ \gen0 \ = \ \big\{\sum n_aa \, \big| \, n_a\in\N\text{ and }n_a=0\text{ for all but finitely many }a\in A-\{0\} \big\}.
\]
An {\em ideal} of $A^+$ is a submodule of $A^+$, i.e.\ a subset $I$ that contains $0$, is closed under addition and that contains $xy$ for all $x\in A^+$ and $y\in I$. We write $I=\gen S$ for the ideal generated by a subset $S$ of $A^+$.

\begin{df}
 A \emph{band} is a pointed monoid $B$ together with a semiring ideal $N_B$ (the \emph{null set}) of $B^+$ such that for every $a\in B$ there is a unique $b\in A$ with $a+b\in N_B$. A {\em band morphism} is a multiplicative map $f:B_1\to B_2$ with $f(0)=0$ and $f(1)=1$ such that $\sum n_af(a)\in N_{B_2}$ for every $\sum n_a a\in N_{B_1}$. This defines the category $\Bands$.
\end{df}

We denote the unique $b$ for which $a+b\in N_B$ by $-a$ and call it the \emph{additive inverse of $a$}. We write $a-b$ for $a+(-b)\in B^+$. The axioms imply that $-0=0$,\quad $(-1)^2=1$ and $-a=(-1)\cdot a$ for all $a\in A$.

The category of bands contains all limits and colimits. Its terminal object is the \emph{zero band} $B=\{0\}$ with $0=1$ and $N_B=\{0\}=B^+$. Its initial object is $\Funpm=\{0,1,-1\}$ with 
\[
 N_\Funpm \ = \ \{n.1+n.(-1)\mid n\geq0\},
\]
which can be seen as a quadratic field extension of $\Fun$.

\subsubsection*{Examples} 
A (commutative and unital) ring $R$ is naturally a band with pointed monoid $B=R$ and null set
\[ \textstyle
 N_B \ = \ \big\{ \sum n_a a \in B^+ \, \big| \, \sum n_a a=0\text{ as equality in $R$} \big\}.
\]
Since a map $f:R_1\to R_2$ is a ring homomorphism if and only if it is a band morphism, this establishes a fully faithful embedding $\Rings\to\Bands$. We typically denote the band with the same symbol as the ring, e.g.\ we use $\F_q$ aned $\Z$ for the band associated with the epinonymous ring.

This embedding extends, in fact, to hyperrings. Namely, a hyperring $R$ (with hypersum $\boxplus$) defines the band with pointed monoid $B=R$ and null set
\[ \textstyle
 N_B \ = \ \big\{ \sum a_i \in B^+ \, \big| \, 0\in \bigboxplus a_i \big\}.
\]
The most relevant examples of hyperrings for our text are the \emph{Krasner hyperfield} $\K=\{0,1\}$ with null set $N_\K=\{n\cdot 1\mid n=0\text{ or }n\geq 2\}$, the \emph{tropical hyperfield} $\T=\R_{\geq0}$ with null set
\[ \textstyle
 N_\T \ = \ \big\{ \sum n_a a \, \big| \, n_a=0\text{ for all }a\in\R_{>0}\text{ or }n_b\geq2\text{ for }b=\max\{a\mid n_a\neq0\}\big\}
\]
and the \emph{hyperring $\cO_\T=[0,1]$ of tropical integers} whose null set is $N_{\cO_\T}=N_\T\cap\cO_\T^+$.

Other examples of bands are partial fields, fuzzy rings and idylls. On the other hand, a band $B$ is naturally an ordered blueprint with underlying monoid $B$, ambient semiring $B^+$ and partial order $\gen{0\leq \sum a_i\mid \sum a_i\in N_B}$.

\subsubsection*{Base extension to rings}
There is a natural functor $(-)_\Z^+:\Bands\to\Rings$, which sends a band $B$ to the ring
\[
 B^+_\Z \ = \ B^+\otimes_\N\Z \ = \ \Z[B]/\gen{N_B}
\]
and a band morphism $f:B\to C$ to the ring homomorphism $f^+\colon B^+_\Z\to C^+_\Z$ with
\[ \textstyle
 f^+\big(\big[\sum n_a\cdot a\big]\big) \ = \ \big[\sum n_a\cdot f(a)\big].
\]
This functor is left adjoint to the fully faithful embedding $\Rings\to\Bands$. In particular, $B_\Z^+\simeq R$ if $B$ is the band associated with the ring $R$.

\subsection{Algebras and quotients}
\label{subsection: Algebras and quotients}

Let $k$ be a band. A \emph{$k$-algebra} is a band $B$ together with a band morphism $\alpha_B:k\to B$. A (\emph{$k$-linear}) {\em morphism} of $k$-algebras $B$ and $C$ is a band morphism $f:B\to C$ for which $\alpha_C=f\circ\alpha_B$. This defines the category $\Alg_k$ of $k$-algebras.

\subsubsection*{Tensor products} 
The tensor product $B\otimes_kC$ is a particular instance of a colimit of bands: its monoid consists of the equivalence classes $b\otimes c$ of elements $(b,c)\in B\times C$ for the equivalence relation $\sim$ generated by $(ab,c)\sim(b,ac)$ for $a\in k$, $b\in B$ and $c\in C$. Its null set is the semiring ideal of $(B\otimes_kC)^+$ generated by expressions of the forms $\sum b_i\otimes 1$ and $\sum 1\otimes c_j$ for which $\sum b_i\in N_B$ and $\sum c_j\in N_C$.

\subsubsection*{Free algebras} 
Let $k$ be a band and $\{T_i\}_{i\in I}$ a set. Let $\gen{T_i}_\N$ be the free monoid generated by the $T_i$, which consists of all finite products of the form $\prod T_i^{e_i}$ with $e_i\geq0$. The \emph{free $k$-algebra in $\{T_i\}$} is the monoid
\[ \textstyle
 k[T_i] \ = \ k[T_i\mid i\in I] \ = \ \big\{ c\prod T_i^{e_i} \, \big| \, c\in k,\, e_i\in\N\text{ with all but finitely many }e_i=0 \big\}
\]
with the identification $0\prod T_i^{e_i}=0\prod T_i^0$ and its null set is generated by the expressions $\sum c_i\prod T_i^0$ for which $\sum c_i\in N_k$. The map $c\mapsto c\prod T_i^0$ defines an injection $k\to k[T_i]$ of bands, and we identify $c\in k$ with $c\prod T_i^0\in k[T_i]$. Note that every map $\{T_i\}\to B$ into a $k$-algebra $B$ extends uniquely to a $k$-linear band morphism $k[T_i]\to B$.

\subsubsection*{Quotients} 
Let $B$ be a band and $S$ a subset of $B^+$. The \emph{quotient of $B$ by $S$} is the band whose monoid 
\[
 \bgenquot BS \ = \ B/\sim
\]
consists of the equivalence classes $\bar a$ of $a\in B$, where $\sim$ is generated by the relations $ac\sim bc$ for which $a-b\in S$ and $c\in B$, and whose null set is 
\[\textstyle
 N_{\bgenquot BS} \ = \ \big\langle \sum n_{\bar a}\bar a \, \big| \, \sum n_aa\in N_B\cup S \big\rangle.
\]
The quotient map $\pi:B\to \bgenquot BS$ with $\pi(a)=\bar a$ is a surjective band morphism, and every band morphism $f:B\to C$ with $\sum n_af(a)\in N_C$ for all $\sum n_aa\in S$ factors uniquely through $\pi$.

\subsubsection*{Examples}
The constructions of free algebras and quotients allow us to write every band in the form $\bgenquot{\Funpm[T_i]}S$. For instance,
\[
 \F_2 \ = \ \bgenquot{\Funpm}{1+1}, \qquad \F_3 \ = \ \bgenquot{\Funpm}{1+1+1}, \qquad \K \ = \ \bgenquot{\Funpm}{1+1,\, 1+1+1},
\]
\[
 \F_4 \ = \ \bgenquot{\Funpm[T]}{1+1,\ T^3+1,\ T^2+T+1}, \qquad \F_5 \ = \ \bgenquot{\Funpm[T]}{T^2+1,\ T+T+1}.
\]

\section{Band schemes}
\label{section: band schemes}

For the versed reader in $\Fun$-geometry, we define a band scheme as an ordered blue scheme $X$ for which $\cO_X(U)$ is a band for all open subsets $U$ of $X$. We spare the reader who is not acquainted with ordered blue schemes the lengthy definitions, but rather work with explicit example classes.

\subsection{Band schemes as functor of points}

For our purposes, we define a \emph{band scheme} as a functor $\cX:\Bands\to\Sets$. This functor is assumed to have an \emph{affine open covering}, but we grossly neglect this condition in this text. A {\em morphism} of band schemes is a morphism of functors. 

An \emph{affine band scheme} is a \emph{representable} functor, i.e.\ a functor isomorphic to $\Hom(B,-)$ for some band $B$. We write $\Spec B=\Hom(B,-)$. By the Yoneda lemma, a { morphism} $\varphi:\Spec C\to\Spec B$ of affine band schemes is induced by a band morphism $f:B\to C$.

Let $k$ be a band. A \emph{$k$-scheme} is a functor $\cX:\Alg_k\to\Sets$. We denote the category of $k$-schemes by $\Sch_k$. Since $\Funpm$ is initial in $\Bands$, a band scheme is the same as an $\Funpm$-scheme, and $\Sch_\Funpm$ corresponds to the category of band schemes.

\subsection{Projective band schemes}

Let $k$ be a band and $n\geq0$. The \emph{projective $n$-space over $k$} is the functor $\P^n_k:\Alg_k\to\Sets$ that sends a $k$-algebra $B$ to the set
\[
 \P^n_k(B) \ = \ \big\{ [x_0:\dotsc:x_n] \, \big| \, (x_0,\dotsc,x_n)\in B^{n+1}\text{ with }x_i\in B^\times\text{ for some }i \big\}
\]
of equivalence classes of $B^{n+1}/B^\times$. Note that a $k$-linear morphism $f:B\to C$ induces the map $\varphi:\P^n_k(B)\to\P^n_k(C)$ with $\varphi([x_0:\dotsc:x_n])=[f(x_0):\dotsc :f(x_n)]$, which defines $\P^n_k$ as a functor.

A formal expression $P=\sum c_{e} T_0^{e_0}\dotsb T_n^{e_n}$ with $c_{e}\in k^+$ and multi-index $e=(e_0,\dotsc,e_n)\in\N^{n+1}$ is called a \emph{homogeneous polynomial over $k$} if there is a $d$ (the \emph{degree}) such that $\sum_{i=0}^n e_i=d$ whenever $c_e\neq0$. Substituting $T_i$ by $x_i$ evaluates $P$ in a tuple $(x_0,\dotsc,x_n)\in k^{n+1}$, and the condition $\sum c_{e} x_0^{e_0}\dotsb x_n^{e_n}\in N_k$ is well defined for equivalence classes $[x_0:\dotsc:x_n]$ if $P$ is homogeneous.

A \emph{projective $k$-scheme} is a subfunctor $\cX$ of $\P^n_k$ of the form
\[ \textstyle
 \cX(B) \ = \ \big\{ [x_0:\dotsc:x_n] \in \P^n_k(B) \, \big| \, \sum \alpha_B(c_{i,e}) x_0^{e_0}\dotsb x_n^{e_n} \, \in \, N_B \text{ for all }i\in I\big\}
\]
where $\alpha_B:k\to B$ is the structure map of the $k$-algebra $B$ and where $\{\sum c_{i,e} x_0^{e_0}\dotsb x_n^{e_n}\}_{i\in I}$ is a set of homogeneous polynomials over $k$. 

In the following, all projective band schemes are defined over $k=\Funpm$, which allows us to omit the subscript $k$. In particular, we write $\P^n$ for $\P^n_\Funpm$.

\subsubsection*{Grassmannians}

Let $E=\{1,\dotsc,n\}$ and $0\leq r\leq n$. Let $\binom Er$ be the family of $r$-subsets of $E$. The Grassmannian over $\Funpm$ is the functor $\Gr(r,n):\Bands\to\Sets$ that sends a band $B$ to the set $\Gr(r,n)(B)$ of equivalence classes $[x_I]_{I\in\binom Er}\in B^N/B^\times$ for which $x_I\in B^\times$ for some $I\in\binom Er$ and that satisfy the \emph{Pl\"ucker relation}
\[
 \sum_{k=0}^r \, (-1)^{\epsilon(j_k,I)+\epsilon(j_k,J)} \, \cdot \, x_{J-\{j_k\}} \, \cdot \, x_{J'\cup\{j_k\}} \ \ \in \ \ N_B 
\]
for all $J=\{j_0,\dotsc,j_{r}\}$ and $J'=\{j'_1,\dotsc,j'_{r-1}\}$ where $\epsilon(j,I)=\#\{i\in I\mid i<j\}$ and $x_{I}=0$ if $\# I< r$. Choosing a linear order on $\binom Er$ defines an inclusion $\Gr(r,n)$ into $\P^N$ for $N=\binom nr-1$, which we call the \emph{Pl\"ucker embedding} and which gives $\Gr(r,n)$ the structure of a projective band scheme. 

If $k$ is a field, then $\Gr(r,n)(k)$ stays naturally in bijection with the usual Grassmannian. 
We recover other well-known spaces: 
\begin{enumerate}
 \item $\Gr(r,n)(\K)$ is canonically bijective to the set of all matroids of rank $r$ on $E=\{1,\dotsc,n\}$; 
 \item $\Gr(r,n)(\T)$ is canonically bijective to the Dressian $\Dr(r,n)$;
 \item $\Gr(r,n)(F)$ is canonically bijective to all $F$-matroids of rank $r$ on $E$ for every idyll $F$.
\end{enumerate}
For more details, we refer the reader to \cite{Baker-Lorscheid21}.

\subsubsection*{Flag varieties}

Also flag varieties extend to band schemes; cf.\ \cite{Jarra-Lorscheid22} for more details. Given a finite set $E=\{1,\dotsc,n\}$ and a tuple $\br=(r_1,\dotsc,r_s)$ of increasing integers $r_1,\dotsc,r_s\in E$ (the \emph{type}), we define the flag variety $\Fl(\br,E)=\Fl(r_1,\dotsc,r_s,E)$ as the functor $\Bands\to\Sets$ that sends a band $B$ to the set of all tuples 
\[
 \big([x_{r_1,I}],\dotsc,[x_{r_s,I}]\big) \quad \in \quad \prod_{i=1}^s \Gr(r_i,E)(B)
\]
that satisfy the incidence relations
\[
 \sum_{k=0}^{r_{i'}} \, (-1)^{\epsilon(j_k,J)+\epsilon(j_k,J')} \, \cdot \, x_{r_{i'},J-\{j_k\}} \, \cdot \, x_{r_i,J'\cup\{j_k\}} \ \ \in \ \ N_B 
\]
for all $i\leq i'$ and $J=\{j_0,\dotsc,j_{r_{i'}}\}$ and $J'=\{j'_2,\dotsc,j'_{r_i}\}$ where we define $x_{r,I}=0$ if $\# I< r$.

\section{Crowds}
\label{section: crowds}

In the existing approaches to algebraic groups over $\Fun$, the group law $\mu:G\times G\to G$ of an algebraic group $G$ does not descend without some yoga to allow for such morphisms over $\Fun$; cf.\ the solutions in \cite{Lorscheid12b} and \cite{Lorscheid18}. The reason for this difficulty is that matrix groups typically involve the additive structure of the tensor powers of its coordinate algebra. In contrast to usual algebraic geometry where tensor products are rings, tensor products fail to be closed under addition in $\Fun$-geometry.

We present in this text a different viewpoint that aligns with the very idea of bands (and, more generally, ordered blueprints): instead of considering the group law as a function $\mu:G\times G\to G$, we consider it as a collection of triples $(a,b,c)$ whose product is $1$, i.e.\ $abc=1$, respectively. This allows us to recover the product $ab=\mu(a,b)$ as the unique element $c$ for which there is a $d$ with $1cd=1$ and $abd=1$. We formalize this type of structure in the following.

\subsection{Definitions}
\label{subsection: definition of crowds}

\begin{df}\label{def: crowd}
 A \emph{crowd} is a set $G$ together with a element $1\in G$ (the \emph{identity}) and a subset $R\subset G^3$ (the \emph{crowd law}) that satisfy the following axioms for all $a,b,c\in G$:
 \begin{enumerate}[label = (C\arabic*)]\setcounter{enumi}{0}
  \item\label{C1} $(a,1,1)\in R$ if and only if $a=1$;
  \item\label{C2} $(a,b,1)\in R$ implies $(b,a,1)\in R$;
  \item\label{C3} $(a,b,c)\in R$ implies $(c,a,b)\in R$.
 \end{enumerate}
 A \emph{crowd morphism} is a map $f:G_1\to G_2$ between crowds $G_1$ and $G_2$ such that $f(1)=1$ and $\big(f(a),f(b),f(c)\big)\in R_2$ for all $(a,b,c)\in R_1$. This defines the category $\Crowds$.
 
 Let $G$ be a crowd. The \emph{inversion of $G$} is the subset
 \[
  R^{(2)} \ = \ \big\{ (a,b)\in G\times G \, \big| \, (a,b,1)\in R \big\}
 \]
 of $G^2$. For $a,b\in G$, the \emph{inverse of $a$} is the subset
 \[
  a^{-1} \ = \ \big\{ b\in G \, \big| \, (a,b)\in R^{(2)} \big\} 
 \]
 of $G$ and the \emph{product of $a$ and $b$} is the subset
 \[
  ab \ = \ \big\{ c\in G \, \big| \, (a,b,d)\in R\text{ and }(d,c)\in R^{(2)}\text{ for some }d\in G \big\}.
 \]
 of $G$.
\end{df}
 
\begin{rem}\label{rem: first properties of crowds}
 We can formulate the axioms of a crowd in a more symmetric form as follows. Let $G$ be a set, $1\in G$ and $R\subset G^3$. Then $G$, together with $1$ and $R$, is a crowd if and only if it satisfies the following axioms.
 \begin{enumerate}[label = (C\arabic*${}^\ast$)]\setcounter{enumi}{0}
  \item The following conditions for $a\in G$ are equivalent:
  \[
   (a) \quad a=1; \qquad (b) \quad (a,1,1)\in R; \qquad (c) \quad (1,a,1)\in R; \qquad (d) \quad (1,1,a)\in R.
  \]
  \item The following conditions for $a,b\in G$ are equivalent:
  \[   
   \begin{aligned}
    &(a) \quad (a,b,1)\in R; \qquad &&(b) \quad (b,1,a)\in R; \qquad &&(c) \quad (1,a,b)\in R; \\
    &(d) \quad (b,a,1)\in R; \qquad &&(e) \quad (a,1,b)\in R; \qquad &&(f) \quad (1,b,a)\in R.
   \end{aligned}
  \]
  \item The following conditions for $a,b,c\in G$ are equivalent:
  \[
   (a) \quad (a,b,c)\in R; \qquad (b) \quad (c,a,b)\in R; \qquad (c) \quad (b,c,a)\in R.
  \]
 \end{enumerate}
\end{rem}

\begin{rem}
 Additional properties that are satisfied by many crowds of interest are the following:
 \begin{enumerate}[label = (C\arabic*)]\setcounter{enumi}{3}
  \item\label{C4} $a^{-1}\neq\emptyset$ for all $a\in G$.
  \item\label{C5} $(a,b,c),\ (a,\bar a,1),\ (b,\bar b,1),\ (c,\bar c,1)\in R$ implies $(\bar c,\bar b,\bar a)\in R$ for all $a,\bar a,b,\bar b,c,\bar c\in G$.
 \end{enumerate}
 The concept of abelian groups extends naturally to crowds in the form of the axiom:
 \begin{enumerate}[label = (C\arabic*)]\setcounter{enumi}{5}
  \item\label{C6} $(a,b,c)\in R$ implies $(b,a,c)\in R$ for all $a,b,c\in G$.
 \end{enumerate}
\end{rem}

We collect some additional properties of a crowd without proof.

\begin{lemma}
 Let $G$ be a crowd. Then: 
 \begin{itemize}
  \item[{\rm (1)}] $1^{-1}=\{1\}$ and $1\notin a^{-1}$ for $a\neq1$;
  \item[{\rm (2)}] $1\in ab$ if and only if $b\in a^{-1}$.
 \end{itemize}
\end{lemma}

\begin{rem}
 The category of crowds has good properties such as the following. The \emph{trivial crowd} $G=\{1\}$ with $R=\{(1,1,1)\}$ is the initial and terminal object in $\Crowds$. The functor $\Crowds\to\Sets_\ast$ to pointed sets that sends a crowd $(G,1,T)$ to the pointed set $G$ with base point $1$, has a fully faithful left adjoint $\Sets_\ast\to\Crowds$, which sends a pointed set $X$ with base point $\ast$ to the \emph{free crowd $G[X]=(X,\ast,T)$ on $X$} with $T=\{(1,1,1)\}$.
\end{rem}

\begin{rem}
 There is a natural alternative to axiomatize a crowd: consider a set $G$ and subsets $R^{(i)}\subset G^i$ for $i=1,2,3$ such that $R^{(1)}=\{1\}$ is a singleton. Consider $G^i$ embedded in $G^{i+1}$ as $G^i\times\{1\}$. Then the crowd axioms \ref{C1}--\ref{C3} are equivalent to 
 \begin{enumerate}
  \item $R^{(i)}$ is invariant under cyclic permutation of the factors for $i=1,2,3$;
  \item $R^{(i)}=R^{(i+1)}\cap G^i$ as subsets of $G^{i+1}$ for $i=1,2$.
 \end{enumerate}
 This reformulation of the crowd axioms suggests a generalization of crowds to relations $R^{(i)}\subset G^i$ for larger $i$ (e.g.\ all $i\in\N$). At the moment of writing, we do not see the need for such a generalization and dismiss this viewpoint from this text.
\end{rem}

\subsection{Groups as crowds}
\label{subsection: groups as crowds}

The principal example of crowds are groups: given a group $G$ with neutral element $1$, we define $R=\{(a,b,c)\in G^3\mid abc=1\}$. Then $G$ together with $1$ and $R$ is a crowd: \ref{C1} is evident; $(a,b,1)\in R$ if and only if $b=a^{-1}$ is the (unique) inverse of $a$, and thus also $(b,a,1)\in R$, which shows \ref{C2}; $(a,b,c)\in R$ if and only if $c^{-1}=ab$ is the (unique) product of $a$ and $b$, and thus also $(c,a,b)\in R$, which shows \ref{C3}.

This discussion should also serve as some intuition for the definition of $a^{-1}$ and $ab$ for more general crowds: for groups, these sets are singletons, and their respective elements coincide with the usual inverse and product of elements of a group.

We sharpen these insights into the relation between groups and crowds in the following result.

\begin{proposition}\label{prop: groups as crowds}
 A map $f:G_1\to G_2$ between groups is a group homomorphism if and only if it is a crowd morphism. This yields a fully faithful embedding $\Groups\to\Crowds$. A crowd $G$ is isomorphic to a group if and only if it satisfies the property that for all $a,b,c,d,e\in G$ we have:  
 \begin{enumerate}
  \item[{\rm (1)}] \label{group1} $a^{-1}$ and $ab$ are singletons;
  \item[{\rm (2)}]\label{group2} $ad=ec$ if $d\in bc$ and $e\in ab$.
 \end{enumerate}
\end{proposition}

\proof
 We begin with the second claim. Let $G$ be a group and $R$ the associated crowd law. We refer to the inverse element of $a\in G$ by $\iota(a)$ and to the product of $a$ and $b$ by $\mu(a,b)$ to differentiate it from the subsets $a^{-1}=\{\iota(a)\}$ and $ab=\{\mu(a,b)\}$ that belong to the crowd $(G,1,R)$. This verifies, in particular, property \eqref{group1}. In order to verify property \eqref{group2}, consider $d\in bc$ and $e\in ab$. Then $d^{-1}=\{\iota(d)\}$ and $e^{-1}=\{\iota(e)\}$ as well as $\mu(\mu(b,c),\iota(d))=1$ and $\mu(\mu(a,b),\iota(e))=1$. Thus $d=\mu(b,c)$ and $e=\mu(a,b)$. By \eqref{group1},
 \[
  ad \ = \ \{ \mu(a,d) \} \ = \ \{ \mu(a,\mu(b,c))\} \ = \ \{ \mu(\mu(a,b),c)\} \ = \ \{\mu(e,c)\} \ = \ ec
 \]
 where the middle equality follows from the associativity of $\mu$. Thus \eqref{group2}.
 
 Conversely, assume that a crowd $(G,1,R)$ satisfies \eqref{group1} and \eqref{group2}. We define $\mu(a,b)$ as the unique element of $ab$ and $\iota(a)$ as the unique element of $a^{-1}$. We claim that $(G,\mu,1)$ is a group with inversion $\iota$. 
 
 We verify the group axioms, starting with the associativity of $\mu$. Consider $a,b,c\in G$. Let $d=\mu(b,c)$ and $e=\mu(a,b)$, i.e.\ $d\in bc$ and $e\in ab$. By \eqref{group2},
 \[
  \{ \mu(a,\mu(b,c))\} \ = \ \{ \mu(a,d) \} \ = \ ad \ = \ ec \ = \ \{\mu(e,c)\} \ = \ \{ \mu(\mu(a,b),c)\},
 \]
 which implies the associativity of $\mu$ since the sets in this equality are singletons. In order to prove that $1$ is right neutral, we consider $a\in G$ and $a1=\{c\mid (a,1,\iota(c))\in R\}$. Then also $(a,\iota(c),1)\in R$ by axioms \ref{C2} and \ref{C3} of a crowd. Thus $\iota(c)\in a^{-1}=\{\iota(a)\}$ and in conclusion $c=a$. This shows that $a1=\{a\}$ and therefore $\mu(a,1)=a$, as desired. To prove that $\iota(a)$ is right-inverse to $a$, note that $(a,\iota(a),1)\in R$ and $\iota(1)=1$ by the definition of $\iota$. Therefore $1\in a\iota(a)=\{\mu(a,\iota(a))\}$, which shows that $\mu(a,\iota(a))=1$. Thus $(G,\mu,1)$ is a group. Note that this group structure for $G$ induces indeed the crowd law $R$ since $(a,b,c)\in R$ if and only if $\mu(a,b)=\iota(c)$ and thus $\mu(\mu(a,b),c)=1$. This concludes the proof of the second claim of the proposition. 
 
 We continue with the first claim. Consider a group homomorphism $f:G_1\to G_2$ and $(a,b,c)\in R_1$, i.e.\ $\mu(\mu(a,b),c)=1$ in $G_1$. Then $\mu(\mu(f(a),f(b)),f(c))=1$ in $G_2$ and thus $(f(a),f(b),f(c))\in R_2$, which shows that $f:G_1\to G_2$ is a crowd morphism. 
 
 Conversely, assume that $f:G_1\to G_2$ is a crowd morphism between groups, i.e.\ the crowd structures of $G_1$ and $G_2$ satisfy \eqref{group1} and \eqref{group2}. Consider $a,b\in G_1$ and $c=\mu(a,b)$, i.e.\ 
 $ab=\{c\}$ (by \eqref{group1}) and $(a,b,\iota(c))\in R_1$. Thus $(f(a),f(b),f(\iota(c))\in R_2$. Since a crowd morphism preserves inverses, i.e.\ $f(\iota(c))=\iota(f(c))$, we find that $\iota(f(c))\in f(a)f(b)$ and thus $\mu(f(a),f(b))=f(c)$, which shows that $f$ is a group homomorphism. In particular, this implies that the functor $\Groups\to\Crowds$ is fully faithful, which completes the proof.
\eop \\

\begin{ex}[Special linear group]\label{ex: special linear group as crowd}
 The leading example in this text is the special linear group $G=\SL_n(R)$ over a ring $R$. Its elements are all $(n\times n)$-matrices $a=(a_{i,j})_{i,j=1,\dotsc,n}$ with coefficients $a_{i,j}\in R$ whose determinant
 \[
  \det(a) \ = \ \sum_{\sigma\in S_n} \ \sign\sigma \, \cdot \ \prod_{k=1}^n \ a_{k,\sigma(k)}
 \]
 is equal to $1$. Its neutral element is the identity matrix $\1=(\delta_{i,j})$ where $\delta_{i,j}$ is the Kronecker symbol. The crowd law $R\subset G^3$ consists of all triples $(a,b,c)$ of elements of $\SL_n(R)$ with product $\1$, i.e.\ which satisfy 
 \[
  \sum_{k,l=1}^n a_{i,k} \cdot b_{k,l} \cdot c_{l,j} \ = \ \delta_{i,j} 
 \]
 for all $i,j=1,\dotsc,n$.
\end{ex}

\subsection{Algebraic crowds}
\label{subsection: algebraic crowds}

We denote by $\cF_G$ the functor $\cF_G:\Crowds\to\Sets$ that sends a crowd $(G,1,R)$ to its underlying set $G$ and by $\cF_R:\Crowds\to\Sets$ the functor that sends $(G,1,R)$ to its crowd law $R$.

\begin{df}
 An \emph{algebraic crowd} is a functor $\cG:\Bands\to\Crowds$ such that both $\cF_G\circ\cG:\Bands\to\Sets$ and $\cF_R\circ\cG:\Bands\to\Sets$ are band schemes. An {algebraic crowd} $\cG$ is \emph{affine} if both $\cF_G\circ\cG$ and $\cF_R\circ\cG$ are affine band schemes.
\end{df}

\begin{rem}\label{rem: algebraic crowds}
 By a standard application of Yoneda's lemma, algebraic crowds can equivalently be characterized as a band scheme $\cX$ (which is $\cF_G\circ\cG$) together with a subscheme $\cR$ of $\cX$ (which is $\cF_R\circ\cG$ as subfunctor of $\cF_G\circ\cG$) and a morphism $\Spec \Funpm\to \cX$ (which corresponds to the unit $1\in \cG(B)$ for all bands $B$) that satisfy certain axioms that correspond to \ref{C1}--\ref{C3}.
 
 If $\cG$ is affine, say $\cF_G\circ\cG\simeq\Hom(C,-)$ and $\cF_R\circ\cG\simeq\Hom(Q,-)$ for bands $C$ (the \emph{coordinate band}) and $Q$, then the crowd structure of $\cG$ is equivalent to epimorphisms $\eta:C\to \Funpm$ (the \emph{counit}) and $\br:C^{\otimes3}\to Q$ (the \emph{colaw}) that satisfy the following axioms for all morphisms $\alpha,\beta,\gamma:C\to B$ into a band $B$ and $\eta_B:C\stackrel\eta\to\Funpm\to B$:
 \begin{enumerate}[label = (C\arabic*${}^*$)]\setcounter{enumi}{0}
  \item The morphism $\alpha\otimes\eta_B\otimes\eta_B:C^{\otimes3}\to B$ factors through $\br:C^{\otimes3}\to Q$ if and only if $\alpha=\eta_B$.
  \item If $\alpha\otimes\beta\otimes\eta_B$ factors through $\br$, then also $\beta\otimes\alpha\otimes\eta_B$ factors through $\br$.
  \item If $\alpha\otimes\beta\otimes\gamma$ factors through $\br$, then also $\gamma\otimes\alpha\otimes\beta$ factors through $\br$.
 \end{enumerate}
\end{rem}

\subsection{The special linear group as an algebraic crowd}
\label{subsection: The special linear group as an algebraic crowd}

The special linear algebraic group $\SL_n$, considered as a functor $\Rings\to\Groups\to\Crowds$, extends to an affine algebraic crowd $\SL_n:\Bands\to\Crowds$ over $\Funpm$, which can be described explicitly as follows, using the notation from Remark \ref{rem: algebraic crowds}.
 
The coordinate band of $\SL_n$ is 
\[
 C \ = \ \bgenquot{\Funpm[T_{i,j}\mid i,j=1,\dotsc,n]}{\det(T_{i,j})-1}.
\]
Its counit is the band morphism $\eta:C\to \Funpm$ that sends $T_{i,j}$ to $\delta_{i,j}$. Its colaw is the projection map $\br:C^{\otimes3}\to Q$ where
\[\textstyle
 Q \ = \ \bgenquot{\Funpm[T^{(m)}_{i,j}\mid m=1,2,3,\, i,j=1,\dotsc,n]}S
\]
with
\[\textstyle
 S \ = \ \big\{ \sum_{k,l}\ T^{(\sigma(1))}_{i,k}T^{(\sigma(2))}_{k,l}T^{(\sigma(3))}_{l,j}-\delta_{i,j} \, \big| \, \sigma\in A_3,i,j=1,\dotsc,n \big\}.
\]
where $A_3\subset S_3$ is the alternating group. It follows from our discussion in Example \ref{ex: special linear group as crowd} that for rings $R$, the crowd $\SL_n(R)$ is the crowd associated with the special linear group of $(n\times n)$-matrices with coefficients in $R$. 

This description of $\SL_n(R)$ extends from rings $R$ to bands $B$. Namely, $\SL_n(B)$ is identified with all $(n\times n)$-matrices $a=(a_{i,j})_{i,j=1,\dotsc,n}$ with coefficients $a_{i,j}\in B$ such that $\det(a_{i,j})-1\in N_B$. The crowd colaw $R(B)$ of $\SL_n(B)$ consists of all triples of matrices $(a^{(1)},a^{(2)},a^{(3)})\in\SL_n(B)^3$ such that
\[
  \sum_{k,l} \ a^{(\sigma(1))}_{i,k}a^{(\sigma(2))}_{k,l}a^{(\sigma(3))}_{l,j}-\delta_{i,j} \in N_B
\]
for all $\sigma\in A_3$ and $i,j=1,\dotsc,n$.

In the light of Proposition \ref{prop: groups as crowds}, this shows that the algebraic crowd $\SL_n:\Bands\to\Crowds$ extends the special linear group from rings to bands. In the following, we describe the crowd $\SL_n(B)$ for some bands $B$ of particular interest.

\begin{ex}\label{F1-points of SL_n}
 As a first example, we consider the regular partial field $B=\Funpm$, which is a subband of the integers $\Z$. This allows us to describe $\SL_n(\Funpm)$ as the subset of matrices in $\SL_n(\Z)$ whose coefficients are in $\Funpm=\{0,\pm 1\}$. A triple of matrices $(a,b,c)$ over $\Funpm$ belongs to the colaw if and only if their product as matrices in $\SL_n(\Z)$ is the identity matrix. This also implies that inverses $a^{-1}$ and products $ab$ are singletons or empty.
 
 As a concrete example, $\SL_2(\Funpm)$ has 20 elements, which are, up to a simultaneous sign change of all coefficients:
 \[
  \matrixtwo1001, \quad
  \matrixtwo1101, \quad
  \matrixtwo1{-1}01, \quad
  \matrixtwo1011, \quad
  \matrixtwo10{-1}1, \quad
  \matrixtwo01{-1}0, \quad
  \matrixtwo11{-1}0, \quad
  \matrixtwo{-1}1{-1}0, \quad
  \matrixtwo01{-1}1, \quad
  \matrixtwo01{-1}{-1}.
 \]
 In this case, all inverses $\matrixtwo abcd^{-1}=\big\{\matrixtwo{d}{-b}{-c}{a}\big\}$ are singletons, but certain products are empty, e.g.\ $\matrixtwo1101\cdot\matrixtwo1101=\emptyset$.
\end{ex}

\begin{ex}\label{K-points of SL_n}
 As the second example, we consider the Krasner hyperfield $\K$. The elements of $\SL_n(\K)$ are $(n\times n)$-matrices $a$ with coefficients $a_{i,j}\in\{0,1\}$ that satisfy the condition $\det(a_{i,j})-1\in N_\K$, which means that $\prod_{i=1}^n a_{i,\sigma(i)}=1$ for some permutation $\sigma\in S_n$. For $n=2$, we have
 \[
  \SL_2(\K) \ = \ \big\{ \
  \matrixtwo1001, \ 
  \matrixtwo0110, \ 
  \matrixtwo1101, \ 
  \matrixtwo1011, \ 
  \matrixtwo1110, \ 
  \matrixtwo0111, \ 
  \matrixtwo1111 \ 
  \big\}.
 \]
 In this case, inverses $\matrixtwo abcd^{-1}=\big\{\matrixtwo dbca\big\}$ are singletons and the product of two matrices $a,b\in\SL_2(\K)$ is given by
 \[\textstyle
  ab \ = \ \big\{ c\in\SL_2(\K) \,\big|\, c_{i,j}+\sum_{k=1}^2 a_{i,k}b_{k,j} \in N_\K\text{ for }i,j=1,2 \big\},
 \]
 which turns $\SL_2(\K)$ into a (non-commutative) canonical hypergroup in the sense of Mittas (\cite{Mittas69}).

 This fact is, however, coincidental to $\SL_2(\K)$. For $n=3$, inverses are not singletons anymore. For example,
 \[
  \matrixthree111111001^{-1} \ = \ \Big\{ \ \matrixthree110111001, \quad \matrixthree111110001, \quad \matrixthree111111001 \ \Big\}.
 \]
\end{ex}

\begin{ex}[Tropical special linear group]
 As a third example, we consider $\SL_2$ for the tropical hyperfield $\T$ and the hyperring $\cO_\T$ of tropical integers. In the latter case $\cO_\T$, we have
 \[
  \SL_2(\cO_\T) \ = \ \big\{ \ \matrixtwo{1}{a}{b}{1},\ \matrixtwo a11b \ \big|\, a,b\in\cO_\T \,\big\},
 \]
 with unique inverses
 \[
  \matrixtwo{1}{a}{b}{1}^{-1} \ = \ \big\{ \matrixtwo{1}{a}{b}{1} \big\}, \qquad \matrixtwo a11b ^{-1} \ = \ \big\{ \matrixtwo{b}{1}{1}{a} \big\},
 \]
 and products
 \[\textstyle
  ab \ = \ \big\{ c\in\SL_2(\cO_\T) \,\big|\, c_{i,j}+\sum_{k=1}^2 a_{i,k}b_{k,j} \in N_{\cO_\T}\text{ for }i,j=1,2 \big\}.
 \]
 In particular, $\SL_2(\cO_\T)$ forms a canonical hypergroup, which extends the corresponding result for $\SL_2(\K)$. The tropical points behave less well:
 \[
  \SL_2(\T) \ = \ \big\{ \ \matrixtwo abcd \ \big| \, ad\leq bc=1\text{ or }bc\leq ad=1\text{ or }1\leq ad=bc \, \big\},
 \]
 which contains tropically singular matrices of the form $\matrixtwo aaaa$ with $a\geq 1$, whose inverse
 \[
  \matrixtwo aaaa^{-1} \ = \ \big\{ \ \matrixtwo bcbc \ \big| \, ab,ac\geq 1 \, \big\}
 \]
 is not a singleton. From the outset, this seems pathological. However, we find back the better behaved set of quasi-invertible tropical matrices with permanent $1$, in the sense of Izhakian and Rowen in \cite{Izhakian-Rowen11}, as the subset
 \[
  \big\{ \, a\in\SL_2(\T) \, \big| \, \det(a)\notin N_\T \, \big\},
 \]
 which should be thought of as an open subset of $\SL_2(\T)$ in the algebro-geometric sense. It is not clear if this subset has a natural structure of a tropical scheme, but it might fit into Friedenberg-Mincheva's framework of tropical adic spaces in \cite{Friedenberg-Mincheva22}.
\end{ex}

\subsection{Descending linear algebraic groups to \texorpdfstring{$\Funpm$}{F1}}
\label{subsection: Descending algebraic groups to F1}

The previous formalism extends to closed subgroups of $\SL_{n,\Z}$, i.e.\ an $\Funpm$-model of a linear algebraic group is defined in terms of a representation on $\A^n$ with determinant $1$.

Namely, a closed immersion $\iota:G\to\SL_{n,\Z}$ of a linear algebraic group $G$ as a subgroup of $\SL_{n,\Z}$ identifies the coordinate algebra of $G$ with a quotient of $\Z[T_{i,j}\mid i,j=1,\dotsc,n]$ by an ideal $I$ that contains $\det(T_{i,j})-1$. We define the affine algebraic crowd $\cG:\Bands\to\Crowds$ as follows: its coordinate algebra is the band
\[\textstyle
 B \ = \ \bgenquot{\Funpm[T_{i,j}\mid i,j=1,\dotsc,n]}{\sum n_a a \mid \sum n_a a\in I}.
\]
Its counit is the morphism $B\to\Funpm$ with $T_{i,j}\mapsto\delta_{i,j}$. Its colaw $R$ is the quotient of $B^{\otimes3}$ by the relations
\[
 \sum_{k,l=1}^n a_{i,k} \cdot b_{k,l} \cdot c_{l,j} \ - \ \delta_{i,j} 
\]
for $i,j=1,\dotsc,n$ and all cyclic permutations of $a$, $b$ and $c$ where $a_{i,j}$, $b_{i,j}$ and $c_{i,j}$ are the classes of $T_{i,j}$ in the respective factors of $B^{\otimes3}$. The following result is immediate from the construction.

\begin{theorem}\label{thm: F1-models of linear algebraic groups}
 The band $B$ together with the counit $B\to\Funpm$ and the colaw $B^{\otimes3}\to R$ defines an affine algebraic crowd $\cG:\Bands\to\Crowds$ such that $\cG(R)$ is the crowd associated with the group $G(R)$ for every ring $R$.
\end{theorem}

\section{Crowd activities}
\label{section: crowd activities}

In this section, we introduce the natural extension of group actions from groups to crowds, which allows us to extend the action of the algebraic group $\SL_{n,\Z}$ on the schemes $\Fl(\br,n)_\Z$ to the realm of band schemes.

\subsection{Definitions}
\label{subsection: definition of crowd activities}

\begin{df}
 A \emph{crowd activity $($of $G$ on $X$$)$} is a triple $(G,X,T)$ where $G$ is a crowd, $X$ a set and $T$ a subset of $G\times X\times X$. We define
 \[
  a.x \ = \ \{ y\in X \mid (a,x,y)\in T\}  \qquad \text{and} \qquad a.S \ = \ \{ y\in X \mid (a,x,y)\in T\text{ for some }x\in S\}
 \]
 for $a\in G$, $x\in X$ and $S\subset X$. A \emph{morphism of crowd activities $(G,X,T)$ and $(G',X',T')$} is a pair $(f,g)$ of a crowd morphism $f:G\to G'$ and a map $g:X\to X'$ such that $(f(a),g(x),g(y))\in T'$ for all $(a,x,y)\in T$. This defines the category $\CrowdActivities$ of crowd activities.
\end{df}

Note that the collection of subsets $a.x$ of $X$ for varying $a\in G$ and $x\in X$ determine $T$ as
\[
 T \ = \ \{ (a,x,y)\in G\times X\times X \mid y\in a.x \}.
\]
Note further that a pair $(f,g)$ of a crowd morphism $f:G\to G'$ and a map $g:X\to X'$ is a morphism between crowd activities $(G,X,T)$ and $(G',X',T')$ if and only if $g(a.x)\subset f(a).g(x)$ for all $a\in G$ and $x\in X$.

\begin{rem}
 The crowd activities of primary interest for this text satisfy some properties, such as the following:
 \begin{enumerate}[label = (A\arabic*)]\setcounter{enumi}{0}
  \item\label{A1} $1.x=\{x\}$ for all $x\in X$.
  \item\label{A2} $y\in a.x$ if and only if $x\in b.y$ for all $x,y\in X$ and $a,b\in G$ with $(a,b,1)\in R$.
  \item\label{A3} $1.x\subset a.(b.(c.x))$ for all $x\in X$ and $a,b,c\in G$ with $(a,b,c)\in R$.
 \end{enumerate}
 We consider it a highly interesting task to study crowd activities with these properties, possibly in combination with suitable additional axioms for the underlying crowds, such as \ref{C4} and \ref{C5}. 
 
 We do not request properties \ref{A1}--\ref{A3} as axioms for crowd activities since they fail to hold for crowd activities of a crowd $G$ on itself by left multiplication, which is given by
 \[
  a.b \ = \ \{c\in G\mid \text{there is a }d\in c^{-1}\text{ such that }(a,b,d)\in R\}.
 \]
 for $a,b\in G$. In general, this crowd activity does satisfy neither of \ref{A1}--\ref{A3}. Note that if $a^{-1}\neq\emptyset$ for all $a\in G$, then \ref{A2} is equivalent to \ref{C5}.
\end{rem}

\subsection{Group actions as crowd activities}
\label{section: group actions as crowd activities}

Crowd activities generalize group actions in the following sense. Let $\GroupActions$ be the category of group actions $\theta:G\times X\to X$ of a group $G$ on a set $X$. A {\em morphism} of group actions $\theta:G\times X\to X$ and $\theta':G'\times X'\to X'$ is a pair $(f,g)$ of a group morphism $f:G\to G'$ and a map $g:X\to X'$ such that $g\circ\theta=\theta'\circ(f,g)$, i.e.\ $g:X\to X'$ is $G$-equivariant with respect to the induced $G$-action $(g,x)\mapsto \theta'(f(g),x)$ of $G$ on $X'$.

Let $\theta:G\times X\to X$ be a group action. The \emph{crowd activity associated with $\theta$} is the crowd activity $(G,X,T)$ with
\[
 T \ = \ \big\{ (a,x,y)\in G\times X\times X \, \big| \, \theta(a,x)=y \big\}.
\]

\begin{proposition}\label{prop: group actions as crowd activities}
 Let $G$ and $G'$ be groups.
 \begin{enumerate}
  \item[{\rm (1)}]\label{act1} A crowd activity $(G,X,T)$ is induced by a group action if and only if it satisfies properties {\rm \ref{A1}} and {\rm \ref{A3}} and if $a.x$ is a singleton for all $a\in G$ and $x\in X$.
  \item[{\rm (2)}]\label{act2} Let $\theta:G\times X\to X$ and $\theta':G'\times X'\to X'$ be group actions with associated crowd activities $(G,X,T)$ and $(G', X', T')$, respectively. Then a pair $(f,g)$ of a group homomorphism $f:G\to G'$ and a map $g:X\to X'$ is a morphism of group actions if and only if it is a morphism of crowd activities. In other words, the association $\theta\mapsto(G,X,T)$ defines a fully faithful embedding $\GroupActions\to\CrowdActivities$. 
 \end{enumerate}
\end{proposition}

\proof
 We begin with claim \eqref{act1}. Assume that $(G,X,T)$ is induced by a group action $\theta:G\times X\to X$. Then $a.x=\{\theta(a,x)\}$ is a singleton for all $a\in G$ and $x\in X$. Moreover, $1.x=\{\theta(1,x)\}=\{x\}$ for all $x\in X$, which verifies \ref{A1}. If $(a,b,c)\in R$, then $abc=1$ as elements of the group $G$. Thus 
 \[
  a.(b.(c.x)) = \{\theta(a,\theta(b,\theta(c,x)))\} = \{\theta(a,\theta(bc,x))\} = \{\theta(abc,x)\} = \{\theta(1,x)\} = \{x\} = 1.x
 \]
 for all $x\in X$, which verifies \ref{A3}.
 
 Conversely, assume that $(G,X,T)$ satisfies \ref{A1} and \ref{A3} and that $a.x$ is a singleton for all $a\in G$ and $x\in X$. Define $\theta(a,x)$ as the unique element $y$ in $a.x$. We claim that $\theta:G\times X\to X$ is a group action, i.e.\ $\theta(1,x)=x$ and $\theta(ab,x)=\theta(a,\theta(b,x))$ for all $a,b\in G$ and $x\in X$. By \ref{A1}, $1.x=\{x\}$ for all $x\in X$, and thus $\theta(1,x)=x$, as required. 
 
 As the next step consider $c\in G$ and $x\in X$, and let $c^{-1}$ be the inverse element of $G$, i.e.\ $(c^{-1},c,1)\in R$. Then by \ref{A1} and \ref{A3}, we have 
 \[
  \theta(c^{-1},\theta(c,x)) \ = \ \theta(c^{-1},\theta(c,\theta(1,x))) \ = \ x.
 \]
 Thus for $a,b\in G$ with product $ab=c$ as elements of $G$, i.e.\ $(a,b,c^{-1})\in R$, and for $y=\theta(c,x)$, we have $y=\theta(a,\theta(b,\theta(c^{-1},y)))$ by \ref{A1} and \ref{A3}. Therefore
 \[
  \theta(a,\theta(b,x)) \ = \ \theta(a,\theta(b,\theta(c^{-1},\theta(c,x)))) \ = \ \theta(1,\theta(c,x)) \ = \ \theta(c,x),
 \]
 which completes the verification that $\theta:G\times X\to X$ is a group action. 
 
 This proves claim \eqref{act1} of the proposition. Claim \eqref{act2} follows at once from the observation that $g(a.x)\subset f(a).g(x)$ is equivalent with $g(\theta(a,x))=\theta(f(a),g(x))$ since both $g(a.x)=\{g(\theta(a,x))\}$ and $f(a).g(x)=\{\theta(f(a),g(x))\}$ are singletons.
\eop \\

\subsection{Algebraic crowd activities}
\label{section: algebraic crowd activities}

Let us review the concept of an algebraic group acting on a scheme. The category of group actions comes with two forgetful functors $\cF_G:\GroupActions\to\Groups$ and $\cF_X:\GroupActions\to\Sets$, which send a group action $\theta:G\times X\to X$ to $G$ and to $X$, respectively. An \emph{algebraic group action} is a functor $\Theta:\Rings\to\GroupActions$ such that $\cF_G\circ\Theta:\Rings\to\Groups$ is an algebraic group and such that $\cF_X\circ\Theta:\GroupActions\to\Sets$ is a scheme.

This concept generalizes to crowd activities in the following way. First note that the forgetful functors $\cF_G$ and $\cF_X$ extend to functors $\cF_G:\CrowdActivities\to\Crowds$ and $\cF_X:\CrowdActivities\to\Sets$ in the obvious way. Crowd activities come with a third forgetful functor $\cF_T:\CrowdActivities\to\Sets$, which sends a crowd activity $(G,X,T)$ to $T$. 

\begin{df}
 An \emph{algebraic crowd activity} is a functor $\Theta:\Bands\to\CrowdActivities$ such that $\cG=\cF_G\circ\Theta$ is an algebraic crowd and such that both $\cX=\cF_X\circ\Theta$ and $\cT=\cF_T\circ\Theta$ are band schemes. We say that \emph{$\cG$ acts on $\cX$ via $\Theta$}.
\end{df}

\begin{rem}
 Note that an algebraic group action $\Theta:\Rings\to\GroupActions$ satisfies the analogon of the last requirement of an algebraic crowd activity automatically: $\cF_T$ sends a group action $\theta:G\times X\to X$ to the set $\{(a,x,y)\in G\times X\times X\mid \theta(a,x)=y\}$, which is bijective to $G\times X$ via $(a,x,y)\mapsto(a,x)$. Thus $\cF_T\circ\Theta$ is isomorphic to the scheme $(\cG\circ\Theta)\times(\cX\circ\Theta)$.
\end{rem}

Given an algebraic crowd $\cG$ and a band scheme $\cX$, an algebraic crowd activity $\Theta$ with $\cG\simeq\cF_G\circ\Theta$ and $\cX\simeq\cF_X\circ\Theta$ is determined by a subfunctor $\cT$ of $\cG\times\cX\times\cX$ that is a band scheme. We describe the algebraic crowd activities of $\SL_n$ on projective spaces, Grassmannians and flag varieties in the following.

\subsubsection*{Crowd activities on projective spaces}
The action of the algebraic group $\SL_{n,\Z}$ on $\P^{n-1}_\Z$ extends to an algebraic crowd activity of $\SL_n$ on $\P^{n-1}$ via the subfunctor $\cT$ of $\SL_n\times\P^{n-1}\times\P^{n-1}$ that sends a band $B$ to 
\[\textstyle
 \cT(B) \ = \ \big\{ (a,x,y) \in \SL_n(B)\times\P^{n-1}(B)\times\P^{n-1}(B) \, \big| \, \sum_k a_{i,k}x_k - y_i \, \in \, N_B \text{ for }i=1,\dotsc,n\big\}.
\]

\subsubsection*{Crowd activities on Grassmannians}
The algebraic crowd activity of $\SL_n$ on $\Gr(r,n)$ is as follows. Given a band $B$, a $B$-matrix $a=(a_{i,j})\in\SL_n(B)$ and $r$-subsets $I$ and $J$ of $E=\{1,\dotsc,n\}$, we define the $(I,J)$-minor of $a$ as
\[
 a_{I,J} \ = \ (-1)^{\epsilon(I)+\epsilon(J)}\det(a_{i,j})_{i\in I,j\in J}.
\]
where $\epsilon(I)$ is the minimal number of transposition on $E$ needed to map $I$ to $\{1,\dotsc,r\}$. Then the algebraic group action of $\SL_n$ on the scheme $\Gr(r,n)$ extends to a algebraic crowd activity $\Theta:\Bands\to\CrowdActivities$ by sending a band to 
\[\textstyle
 \cT(B) \ = \ \big\{ (a,x,y) \in \SL_n(B)\times\Gr(r,n)(B)\times\Gr(r,n)(B) \, \big| \, \sum_J a_{I,J}x_J - y_I \, \in \, N_B \text{ for }I\in\binom Er\big\}.
\]

\subsubsection*{Crowd activities on flag varieties}
Eventually the algebraic crowd activity of $\SL_n$ on the Grassmannians $\Gr(r,n)$ for various $r$ extend to an algebraic crowd activity $\cT$ on flag varieties $\Fl(\br,n)$ of type $\br=(r_1,\dotsc,r_s)$, which is a closed subscheme of $\prod \Gr(r_i,n)$. A triple of elements $a\in\SL_n(B)$ and $(x_{1},\dotsc,x_{r}),\, (y_{1},\dotsc,y_{r}) \in\Fl(\br,n)(B)$ is in $\cT(B)$ if and only if
\[
 \sum_J a_{I,J}x_{r_i,J} - y_{r_i,I} \, \in \, N_B
\]
for all $i=1,\dotsc,r$ and all $I\subset E$ of cardinality $r_i$.

\section{Tits's dream revisited}
\label{section: Tits's dream revisted}

In this section we explain how Tits's proposed geometry over the field with one element appears naturally in the geometry of flag varieties and crowd activities over the Krasner hyperfield. More accurately, Tits's vision on geometry over $\Fun$ appears as the outer layer of Borovik-Gelfand-White's combinatorial flag varieties.

\subsection{Combinatorial flag varieties}
\label{section: Combinatorial flag varieties}

Borovik, Gelfand and White define in \cite{Borovik-Gelfand-White01} the \emph{combinatorial flag varieties $\Omega_{S_n}$} as the order complex of all matroids on $E=\{1,\dotsc,n\}$ with respect to the partial order given by matroid quotients. More explicitly, the simplices of $\Omega_{S_n}$ correspond to partial flag matroids $(N_1,\dotsc,N_s)$ on $E$ with $0<\rk N_1<\cdots<\rk N_s<n$. The faces of a flag matroid $(N_1,\dotsc,N_s)$ are all flag matroids $(N_{i_1},\dotsc,N_{i_t})$ with $1\leq i_1<\dotsb<i_t\leq s$. The dimension of $\Omega_{S_n}$ is $n-2$. 

See Figure \ref{fig: Coxeter complex and combinatorial flag variety for S3} for an illustration of $\Omega_{S_3}$. We label its vertices by a matrix that represents the corresponding matroid, which is possible since all matroids on $3$ elements are regular.

\begin{figure}[htb]
 \[
  \begin{tikzpicture}
   \node (origin) {};
   \foreach \a in {1,...,14}{\draw (\a*360/7+0/7: 2cm) node [draw,circle,inner sep=2pt] (l\a) {};}
   \foreach \a in {1,...,7}{\draw (\a*360/7+180/7: 2cm) node [draw,circle,inner sep=1.8pt,fill] (p\a) {};
                              \draw [-] (p\a) -- (l\a);
                              \setcounter{tikz-counter}{\a};
                              \addtocounter{tikz-counter}{1};
                              \draw [-] (p\a) -- (l\arabic{tikz-counter});
                              \addtocounter{tikz-counter}{2};
                              \draw [-] (p\a) -- (l\arabic{tikz-counter});
                             }
   \draw [-] (p4) -- (l6);                          
   \node at      (0:2.8cm) {$\tinymatrix 110001$};                          
   \node at  (360/7:2.8cm) {$\tinymatrix 100001$};                          
   \node at  (720/7:2.8cm) {$\tinymatrix 100010$};                          
   \node at (1080/7:2.8cm) {$\tinymatrix 001010$};                          
   \node at (1440/7:2.8cm) {$\tinymatrix 100011$};                          
   \node at (1800/7:2.8cm) {$\tinymatrix 101010$};                          
   \node at (2160/7:2.8cm) {$\tinymatrix 110011$};                          
   \node at  (180/7:2.8cm) {$\tinyvector 001$};                          
   \node at  (540/7:2.8cm) {$\tinyvector 100$};                          
   \node at  (900/7:2.8cm) {$\tinyvector 010$};                          
   \node at (1260/7:2.8cm) {$\tinyvector 011$};                          
   \node at (1620/7:2.8cm) {$\tinyvector 111$};                          
   \node at (1980/7:2.8cm) {$\tinyvector 101$};                          
   \node at (2340/7:2.8cm) {$\tinyvector 110$};   
   \node at  (360/7:2.8cm) {$\tinymatrix 100001$};                          
   \node at  (720/7:2.8cm) {$\tinymatrix 100010$};                          
   \node at (1080/7:2.8cm) {$\tinymatrix 001010$};                          
   \node at  (180/7:2.8cm) {$\tinyvector 001$};                          
   \node at  (540/7:2.8cm) {$\tinyvector 100$};                          
   \node at  (900/7:2.8cm) {$\tinyvector 010$};                          
  \end{tikzpicture}
 \] 
 \caption{The combinatorial flag variety $\Omega_{S_3}$}
 \label{fig: combinatorial flag variety for S3}
\end{figure}
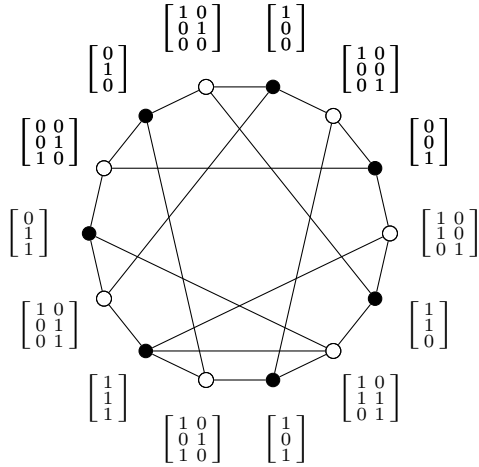

\subsection{Rational points of flag varieties}
\label{section: Rational points of flag varieties}

The band schemes $\Fl(\br,n)$ come with projection maps: if $\br'=(r_{i_1},\dotsc,r_{i_t})$ is a subtype of $\br=(r_1,\dotsc,r_s)$, i.e.\ $1\leq i_1<\dotsb<i_t\leq s$, then the projection $\pi_{\br'}:\Fl(\br,n)\to\Fl(\br',n)$ sends a $B$-rational point $\big([x_{r_1,I}],\dotsc,[x_{r_s,I}]\big)$ to $\big([x_{r_{i_1},I}],\dotsc,[x_{r_{i_t},I}]\big)$.

Given a band $B$, we define a simplicial complex $\Delta_{n}(B)$ as follows. Its simplices are the elements of 
\[
 \coprod_{\br\in\Theta} \ \Fl(\br,n)(B) \quad \text{where} \quad \Theta \ = \ \big\{ \, (r_1,\dotsc,r_s) \, \big| \, s>0,\; 0<r_1<\dotsb<r_s<n \, \big\},
\]
and the dimension of a simplex $\delta\in \Fl(\br,n)(B)$ with $\br=(r_1,\dotsc,r_s)$ is $\dim\delta=s-1$. The faces of $\delta$ are the simplices $\pi_{\br',B}(\delta)$ for the projection $\pi_{\br',B}:\Fl(\br,n)(B)\to\Fl(\br',n)(B)$ for subtypes $\br'$ of $\br$. This defines $\Delta_n$ as a functor from $\Bands$ to simplicial complexes. More to the point, $\Delta_n$ is a simplicial band scheme.

If $B$ is finite, then $\Delta_n(B)$ is finite. In particular, we find:
\begin{itemize}
 \item $\Delta_n(\F_q)$ is the spherical building $\cB_n(\F_q)$ of type $\tA_{n-1}$ over $\F_q$;
 \item $\Delta_n(\K)$ is the combinatorial flag variety $\Omega_{S_{n}}$.
\end{itemize}
The former follows from the well-known identification of $\cB_n$ with flags of linear subspaces of $\F_q^{n}$; the latter follows from the identification of flag matroids of type $\br$ with points of $\Fl(\br,n)(\K)$ (\cite{Jarra-Lorscheid22}).

Due to the functorial definition of $\Delta_n(B)$, every band morphism $f:B\to C$ induces a simplicial map $f_\ast:\Delta_n(B)\to\Delta_n(C)$. In particular, we have a commutative diagram
\[
 \begin{tikzcd}
  \Delta_n(\F_q) \ar[r,"\sim"] \ar[d,"t_{\F_q,\ast}"'] & \cB_n(\F_q) \ar[d,"\mu_q"] \\
  \Delta_n(\K) \ar[r,"\sim"] & \Omega_{S_{n}}
 \end{tikzcd}
\]
of simplicial maps where the horizontal arrows are the aforementioned identifications, $t_{\F_q,\ast}$ is induced by the unique band morphism $t:\F_q\to\K$ with $t(a)=1$ for $a\in\F_q^\times$ and $\mu_q$ sends a flag $(V_1,\dotsc,V_s)$ of subvector spaces of $\F_q^{n+1}$ to the induced flag matroid.

\begin{rem}
 The reader might have noticed the similarity between the combinatorial flag variety $\Omega_{S_3}$ (in Figure \ref{fig: combinatorial flag variety for S3}) with the spherical building of type $\tA_2$ over $\F_2$: except for the $1$-simplex with vertices labelled by $\tinyvector111$ and ${\tinymatrix 110011}$, these two simplicial complexes are equal. This proximity is due to the fact that all flag matroids on $3$ elements but for the mentioned one, are binary. For larger $n$ the discrepancy increases; in fact the percentage of flag matroids that come from a spherical building over any $\F_q$ goes to $0$ as $n$ goes to infinity, by a result of Nelson (\cite{Nelson18}).
\end{rem}

\subsection{Recovering Tits geometries}
\label{section: Recovering Tits geometries}

Borovik, Gelfand and White observe in \cite[section 7.14]{Borovik-Gelfand-White03} that the Coxeter complex $\Pi_{n,\Fun}$ of $S_n$, which Tits envisioned as a geometry over $\Fun$, appears as a subcomplex of the combinatorial flag variety $\Omega_{S_n}=\Delta_{n}(\K)$. We recover this subcomplex from the viewpoint of band schemes as follows: a point $\big([x_{r_1,I}],\dotsc,[x_{r_s,I}]\big)$ of $\Gr(\br,n)(\K)$ lies in $\Pi_{n,\Fun}$ if and only if for $i=1,\dotsc,s$, the tuple $[x_{r_i,I}]$ has precisely one non-zero entry.

More accurately, there is a simplicial subscheme $\Gamma_n$ of $\Delta_n$, which is defined by the vanishing of terms of the form $x_{r,I}x_{r,J}$ with $I\neq J$ and which satisfies the following two properties:
\begin{itemize}
 \item $\Gamma_n(\F_q)$ is the apartment of the canonical basis of $\F_q^{n}$ in $\cB_n(\F_q)$;
 \item $\Gamma_n(\K)$ is the Coxeter complex of $S_{n}$ as a subcomplex of $\Omega_{S_n}$.
\end{itemize}

\subsection{Extending the symmetry group}
\label{section: Extending the symmetry group}

The action of $S_n$ on the Coxeter complex $\Pi_{n,\Fun}=\Gamma_{n}(\K)$ and the combinatorial flag variety $\Omega_{S_n}=\Delta_{n}(\K)$ can be recovered from the crowd activity of $\SL_n$ on $\Fl(\br,n)$ (for various $\br$) as follows.

As a first observation note that the componentwise definition of the crowd activity $\SL_n$ on $\Fl(\br,n)$ implies that the projections $\pi_{\br'}\colon\Fl(\br,n)\to\Fl(\br',n)$ are morphisms of algebraic crowd activities (with respect to the algebraic crowd activities of $\SL_n$ on either flag variety). This defines an algebraic crowd activity of $\SL_n$ on the simplicial band scheme $\Delta_{n}$.

We define $\cN$ as algebraic subcrowd of all monomial matrices of $\SL_n$. As a band scheme, it is defined as
\[
 \cN(B) \ = \ \big\{\, (a_{i,j})\in \SL_n(B) \,\big|\, \text{there is an $\sigma\in S_n$ such that }a_{i,j}=0\text{ whenever }j\neq\sigma(i) \, \big\}
\]
for every band $B$. Its crowd law is the restriction of the crowd law of $\SL_n(B)$ to $\cN(B)$.

Over a ring $R$, the group $\cN(R)$ is the normalizer of the diagonal torus of $\SL_n(R)$. In fact, for every band $B$, the crowd $\cN(B)$ is a group, namely the group of monomial matrices of determinant $1$. In other words, the crowd law of $\cN(B)$ consists of all triples of monomial matrices whose product is $\1$. Since $\K^\times=\{1\}$ and $-1=1$ in $\K$, the group $\cN(\K)$ consists of all permutation matrices and thus $\cN(\K)\simeq S_n$. 

The algebraic crowd activity of $\SL_n$ on $\Delta_n$ restricts to a group action of $\cN$ on $\Delta_n$. Under the identifications $\cN(\K)=\S_n$ and $\Delta_{n}(\K)=\Omega_{S_n}$, this group action recovers the usual action of $S_n$ on $\Omega_{S_n}$.

Turning this observation around, we see that we have extended the group action of $S_n$ on $\Omega_{S_n}$ in a natural way to the crowd activity of $\SL_n(\K)$ on $\Omega_{S_n}$. 

\begin{remark}
{\rm We consider it a highly interesting task to study the properties of this crowd activity and see potential applications to conjectures in combinatorics with a finite field analogon.
}\end{remark}

\subsection{A case study}
\label{section: A case study}

As an explicit example, we study the crowd activity of $\SL_3(\K)$ on $\Omega(S_3)=\Delta_3(\K)$. The elements $a\in \SL_3(\K)$ and $x\in\Delta_3(\K)$ are layered by their ``genericity,'' which heuristically can be thought of as the amount of non-trivial coefficients, and which finds a precise measure in the cardinalities of the orbits $a.x$.

All orbits are non-empty, so the minimal cardinality of an orbit is $1$. For monomial $a_\sigma\in\cN(\K)\subset\SL_3(\K)$ with $a_{\sigma,i,j}=\delta_{i,\sigma(j)}$, all orbits $a_\sigma.x$ are singletons, namely\footnote{Note that we use in the computation for orbits of $a\in\SL_3(\K)$ on elements of $p\in\Gr(2,3)(\K)$ that $a$ acts on matrix representatives $x$ of $p$ in terms of matrix multiplication.}
\[
 a_\sigma.\Big[\begin{smallmatrix} x_1 \\ x_2 \\ x_3 \end{smallmatrix}\Big] \ = \ \Big\{\Big[\begin{smallmatrix} x_{\sigma(1)} \\ x_{\sigma(2)} \\ x_{\sigma(3)} \end{smallmatrix}\Big] \Big\}
 \qquad \text{and} \qquad
 a_\sigma.\Big[\begin{smallmatrix} x_{1,1} & x_{1,2} \\ x_{2,1} & x_{2,2} \\ x_{3,1} & x_{3,2}  \end{smallmatrix}\Big] \ = \ \Big\{\Big[\begin{smallmatrix} x_{\sigma(1),1} & x_{\sigma(1),2} \\ x_{\sigma(2),1} & x_{\sigma(2),2} \\ x_{\sigma(3),1} & x_{\sigma(3),2}  \end{smallmatrix}\Big]\Big\} .
\]
Similarly, for points $x\in\Gamma_3(\K)$, the orbits $a.x$ are singletons for all $a\in\SL_3(\K)$. For instance, 
\[
 a.\Big[\begin{smallmatrix} 1 \\ 0 \\ 0 \end{smallmatrix}\Big] \ = \ \Big\{\Big[\begin{smallmatrix} a_{1,1} \\ a_{2,1} \\ a_{3,1} \end{smallmatrix}\Big] \Big\}
 \qquad \text{and} \qquad
 a.\Big[\begin{smallmatrix} 1 & 0 \\ 0 & 1 \\ 0 & 0 \end{smallmatrix}\Big] \ = \ \Big\{\Big[\begin{smallmatrix} a_{1,1} & a_{1,2} \\ a_{2,1} & a_{2,2}  \\ a_{3,1} & a_{3,2}  \end{smallmatrix}\Big] \Big\}.
\]

These are the only $a\in\SL_3(\K)$ and $x\in\Delta_3(\K)$, respectively, for which all orbits are singletons. Other combinations yield orbits of larger sizes. For instance, $a=\Big[\begin{smallmatrix} 1 & 1 & 0 \\ 0 & 1 & 0 \\ 0 & 0 & 1\end{smallmatrix}\Big]$ has orbits of sizes $1$ and $2$, such as
\[
 \Big[\begin{smallmatrix} 1 & 1 & 0 \\ 0 & 1 & 0 \\ 0 & 0 & 1\end{smallmatrix}\Big].\Big[\begin{smallmatrix} 1 \\ 0 \\ 1\end{smallmatrix}\Big] \ = \ \Big\{ \Big[\begin{smallmatrix} 1 \\ 0 \\ 1\end{smallmatrix}\Big] \Big\}
 \qquad \text{and} \qquad
 \Big[\begin{smallmatrix} 1 & 1 & 0 \\ 0 & 1 & 0 \\ 0 & 0 & 1\end{smallmatrix}\Big].\Big[\begin{smallmatrix} 1 \\ 1 \\ 1\end{smallmatrix}\Big] \ = \ \Big\{ \Big[\begin{smallmatrix} 0 \\ 1 \\ 1\end{smallmatrix}\Big], \, \Big[\begin{smallmatrix} 1 \\ 1 \\ 1\end{smallmatrix}\Big] \Big\}.
\]
The orbits of $a=\Big[\begin{smallmatrix} 1 & 1 & 1 \\ 1 & 1 & 1 \\ 1 & 1 & 1\end{smallmatrix}\Big]$ that are not singletons have cardinality $7$, e.g.\
\[
 \Big[\begin{smallmatrix} 1 & 1 & 1 \\ 1 & 1 & 1 \\ 1 & 1 & 1\end{smallmatrix}\Big]. \Big[\begin{smallmatrix} 1 \\ 1 \\ 1\end{smallmatrix}\Big] \ = \ \Gr(1,3)(\K)
 \qquad \text{and} \qquad
 \Big[\begin{smallmatrix} 1 & 1 & 1 \\ 1 & 1 & 1 \\ 1 & 1 & 1\end{smallmatrix}\Big]. \Big[\begin{smallmatrix} 1 & 0 \\ 0 & 1 \\ 1 & 1\end{smallmatrix}\Big] \ = \ \Gr(2,3)(\K).
\]

\subsection{Other Dynkin types}
\label{section: Other Dynkin types}

In this text, we have demonstrated in the sample case of $\SL_n$ how our proposed solution to Tits's dream leads to interesting geometric structures with the potential for future applications.

We expect that the formalism of this text extends to other $\Fun$-models of $\SL_n$ and other algebraic groups. The subtlety here is that different linear presentations of the same algebraic group lead to different $\Fun$-models. In this sense, {$\Fun$-geometry encodes the representation theory of an algebraic group.} 

In particular, we expect that Tits geometries of Dynkin types $\tB$, $\tC$ and $\tD$ appear as a subcomplex of the space of $\K$-rational points of suitable simplicial band schemes that come equipped with suitable algebraic crowd activities. The space of all $\K$-rational points should reflect the concept of Coxeter matroids by Borovik, Gelfand and White (cf.\ \cite{Borovik-Gelfand-White03}).

We hope that these remarks stimulate a rigorous treatment of Tits geometries, their $\Fun$-models and the relation to matroid theory in the proposed language of band schemes and algebraic crowds.


\ 
\part{\texorpdfstring{$\Fun$}{F1}-polygons}

In the previous sections, we have studied various guises of $\Fun$-geometries under the general umbrella of schemes. In classical scheme theory, say | for the sake of convenience (but without loss of generality) | on the affine level, we consider the category $\ComRing$ of commutative rings with multiplicative identity, which comes with an intitial object $\mathbb{Z}$. As the ring of integers maps uniquely to any commutative ring $A$, applying the controvariant functor $\Spec$, we obtain a diagram 
\begin{equation}
    \Spec(A) \ \epiarrow{1.1cm} \ \Spec(\mathbb{Z}).
\end{equation}

The typical philosophy which one has in mind when dreaming about $\Fun$-geometry, is a relaxation of the category $\ComRing$, in which an object $\Fun$ arises which sits under $\mathbb{Z}$:
\begin{equation}
    \Fun \ \monoarrow{1.1cm}\ \mathbb{Z} \ \monoarrow{1.1cm}\ A, 
\end{equation}
where $A$ is now an object in the new category (and which in particular should be allowed to be a commutative ring), such that applying a new contravariant functor $\widehat{\Spec}$ which appropriately generalizes the functor $\Spec$, gives the desired diagram of projections
\begin{equation}
    \widehat{\Spec}(A) \ \epiarrow{1.1cm} \ \widehat{\Spec}(\mathbb{Z})\ \epiarrow{1.1cm} \ \widehat{\Spec}(\Fun).
\end{equation}

This formalism is exactly what we have studied in the first part of this paper, in the context of the category of bands (instead of commutative rings) and band schemes (instead of classical schemes), and applied to algebraic groups. In this section, we want to consider an incidence-geometrical analogon of the aforementioned formalism. In particular, we want to consider epimorphisms 
\begin{equation}
\label{epidiag}
    \mG \ \epiarrow{1.1cm} \ \Delta,
\end{equation}
where $\mG$ is an object in some category $\mathsf{Geom}$ of combinatorial geometries, and $\Delta$ is an ``$\Fun$-object'' in this category. Of course, it is not always clear what an $\Fun$-object in a (geometric) category really is (and this is certainly a question which should be considered in a future paper). 
The leading example, just like in pretty much most of the paper, is the case where $\mathsf{Geom}$ is the category of projective spaces (that is, buildings of type $\tA_N$) of fixed dimension $N$ in a fixed characteristic, and $\Delta$ is the well-defined thin version of these geometries (a complete graph of size $N + 1$). In Thas and Thas \cite{part3,part4}, the case was handled for dimension $N = 2$ and with the additional assumption that the considered projective geometries are {\em finite}; this condition was also shown to be nessecary since free constructions are possible when infinite dimensions are allowed. Diagrams such as (\ref{epidiag}) endow the geometry $\mG$ with what we call an ``$\Fun$-structure'' (and that is also what their scheme-theoretic cousins do). In Thas and Thas \cite{part3,part4}, not only diagrams (\ref{epidiag}) were considered with source a (possibly non-classical) projective plane | they also considered buildings of type $\tB_2$, $\tC_2$ and $\mathsf{I}_4(2)$. Again, the finiteness conditions were necessary.

In this second part of the paper, we will completely determine the epimorphisms 
\begin{equation}
\label{epidiag2}
   \upepsilon:\ \mathbb{P}^3(\F_q) \ \epiarrow{1.1cm} \ \Delta = \mathbb{P}^3(\Fun) = K(4),
\end{equation}
where $\F_q$ is any finite field, and the target is the complete graph on $4$ vertices.  

We will first have a deeper look at the theory of thin generalized polygons, and show that | through a doubling procedure | the structure of thin polygons (which one should see as $\Fun$-versions of thick polygons) is much more elaborate than what one might be tended to think at first. In particular, as soon as the (even) gonality of the polygon is at least $6$, the structure of thin generalized $n$-gons becomes highly complex, which contrasts the cited passage from \cite{Borovik-Gelfand-White01} in the introduction. As we will see, through the doubling procedure applied to a generalized $n$-gon $\Upgamma$ of order $(s,s)$, we obtain a thin polygon $\Upgamma^\Updelta$ with gonality $2n$, and it will be clear that $\Upgamma^\Updelta$ is actually the flag variety of $\Upgamma$.

\medskip
\section{Generalized polygons}

Let $\Upgamma = (\mP,\mL,\I)$ be a point-line geometry, and let $m$ be a positive integer at least $2$. We say that $\Upgamma$ is a {\em weak generalized $m$-gon} if: 
\begin{itemize}
\item[(A)]
any two elements in $\mP \cup \mL$ are contained in at least one ordinary sub $m$-gon (as a subgeometry of $\Upgamma$), and 
\item[(B)]
if $\Upgamma$ does not contain ordinary sub $k$-gons with $2 \leq k < m$. 
\end{itemize} 

For $m = 2$ every point is incident with every line. If $m \geq 3$, we say $\Upgamma$ is a {\em generalized $m$-gon} if furthermore: 
\begin{itemize} 
\item[(C)]
$\Upgamma$ contains an ordinary sub $(m + 1)$-gon as a subgeometry. 
\end{itemize} 

\begin{remark}{\rm 
Note that the generalized $3$-gons are precisely the (axiomatic) projective planes. Generalized $4$-gons, resp. $6$-gons, resp. $8$-gons are also called {\em generalized quadrangles}, resp. {\em hexagons}, resp. {\em octagons}.} 
\end{remark} 

\subsection{Thick and thin polygons}

Equivalently, a weak generalized $m$-gon with $m \geq 3$ is a generalized $m$-gon if it is {\em thick}, meaning that every point is incident with at least three distinct lines and every line is incident with at least three distinct points. A weak generalized $m$-gon is {\em thin} if it is not thick; in that case, we also speak of {\em thin generalized $m$-gons}. If we do not specify $m$ (the ``gonality''), we speak of {\em (weak) generalized polygons}.

\subsection{Order} 

It can be shown that generalized polygons have an {\em order} $(u,v)$: there exists positive integers $u \geq 2$ and $v \geq 2$ such that each 
point is incident with $v + 1$ lines and each line is incident with $u + 1$ points. We say that a weak generalized polygon is {\em finite} if its number of points and lines is finite | otherwise it is {\em infinite}.  If a thin weak generalized polygon has an order $(1,u)$ or $(u,1)$ it is called a {\em thin 
generalized polygon} of order $(1,u)$ or $(u,1)$. 

\ex
A generalized $m$-gon of order $(1,1)$ is also called an {\em apartment}; it is an ``ordinary $m$-gon'' in property (A) in the definition of weak generalized $m$-gon above.  

By the following remarkable theorem, the gonality of generalized polygons is severely restricted if the geometries are finite.

\begin{theorem}[Feit and Higman \cite{FeHi}]
Let $\Upgamma$ be a finite weak generalized $n$-gon of order $(u,v)$ with $n \geq 3$. Then we have one of the following possibilities:
\begin{itemize}
\item[{\rm (1)}]
$\Upgamma$ is an ordinary $n$-gon (so that $u = v = 1$);
\item[{\rm (2)}]
$\Upgamma$ is thick and $n \in \{ 3, 4, 6, 8 \}$;
\item[{\rm (3)}]
$n = 12$ and either $u = 1$ or $v = 1$.
\end{itemize}
\end{theorem}

\subsection{Morphisms and epimorphisms}

A {\em morphism} from a weak generalized polygon $\Upgamma = (\mP,\mL,\I)$ to a weak generalized polygon $\Upgamma = (\mP',\mL',\I')$ is a map $\upalpha: \mP \cup \mL \mapsto \mP' \cup \mL'$ which maps points to points, lines to lines and which preserves the incidence relation (note that we do not ask the gonalities to be the same). We say that a morphism $\upalpha$ is an {\em epimorphism} if $\upalpha(\mP) = \mP'$ and $\upalpha(\mL) = \mL'$.

If an epimorphism is injective, and if the inverse map is also a morphism, then we call it an {\em isomorphism}. An isomorphism of type $\Upgamma: \Upgamma \mapsto \Upgamma$ with $\Upgamma$ a weak generalized polygon, is called an {\em automorphism} of $\Upgamma$. Note that the set $\Aut(\Upgamma)$ of all automorphisms of $\Upgamma$ naturally forms a group under the composition of maps.  

\begin{remark}{\rm 
In categorical language, an {\em epimorphism} is any morphism which is right-cancellative. 
In the category of sets, this is trivially equivalent to asking that the morphism (map) is surjective. Since morphisms between generalized polygons are defined by the underlying maps between the point sets and line sets, it follows that in the categorical sense, epimorphisms between polygons are indeed as above. }
\end{remark}

\section{Thin polygons and \texorpdfstring{$\Fun$}{F1}-polygons}

Let $m \geq 3$ be a positive integer. 
In Tits's  original approach, the $\Fun$-version of a generalized $m$-gon $\Upgamma$ was isomorphic to any apartment of $\Upgamma$ | in other words, it was a weak generalized $m$-gon of order $(1,1)$.  
Since weak generalized $m$-gons with two points per line also merit to be associated to $\Fun$, we say that a weak generalized  $m$-gon is {\em defined over $\Fun$} if it has order $(1,t)$. We say it is {\em strictly defined over $\Fun$} if $t = 1$. (Note that this terminology differs a bit from that used in \cite{chapcomb}.) 

Note that weak generalized polygons are defined (in the approach above) through a property of forbidden subconfigurations, and a covering property by subpolygons (strictly) defined over $\Fun$ (which are all isomorphic). 

We say that $(\Upgamma,\Updelta,\upepsilon)$ is an {\em $\Fun$-generalized polygon} or {\em $\Fun$-polygon} if $\Upgamma$ is a thick generalized $n$-gon, $\Updelta$ a thin 
generalized $n$-gon (of order $(s,1)$), and 
\begin{equation}
\upepsilon:\ \Upgamma\ \epiarrow{1.1cm}\ \Updelta
\end{equation} 
 a surjective morphism.  
We say that $(\Upgamma,\Updelta,\upepsilon)$ is a {\em proper} $\Fun$-polygon if $\Updelta$ has order $(1,1)$. $\Fun$-Polygons are called {\em finite} if $\Upgamma$ (and hence $\Updelta$, as $\upepsilon$ is surjective) is finite. 
Note that if $n = 3$, $s$ necessarily equals $1$.

\subsection{Doubling and polygons defined over \texorpdfstring{$\Fun$}{F1}}

In this subsection and the next one, we only consider finite generalized $n$-gons with $n \geq 3$ (unless otherwise stated), so that $n \in \{ 3, 4, 6, 8\}$ by the Feit-Higman result. 

Due to the fact that one of the parameters of a thin generalized $n$-gon of order $(1,t)$ is $1$ --- and perhaps also due to the use of the suggestive word ``thin'' --- one might have the impression that these geometries have a simple structure. In this subsection and the next one, we will show that this belief is (totally) wrong: only in the cases $n = 3$ and $n = 4$ we end up with easy-to-understand geometries, but the cases $n = 6$ and $n = 8$ are very difficult. 
We first need to explain the \ul{doubling procedure.} 

Let $\Upgamma = (\mP, \mL, \I)$ be a (not necessarily finite) generalized $n$-gon of order $(s,s)$ (for $n = 3$, projective planes of order $(1,1)$ are allowed). Define the {\em double of $\Upgamma$} as the generalized $2n$-gon $\Upgamma^\Updelta$ which arises by letting its point set to be $\mP \cup \mL$, and letting its line set be the flag set of $\Upgamma$ (the set of incident point-line pairs).  Its parameters are $(1, s)$. Note that in some sense, $\Upgamma^\Updelta$ is the (dual) {\em flag-variety} of $\Upgamma$. 

 The full automorphism group of $\Upgamma^\Updelta$ is isomorphic to the  group consisting of all automorphisms and dualities (anti-automorphisms) of $\Upgamma$. Sometimes we prefer to work in the point-line dual of $\Upgamma^\Updelta$, but we use the same notation (while making it clear in what setting we work). 
Vice versa, if $\Upgamma'$ is a thin generalized $2n$-gon of order $(1,s)$, then it is isomorphic to the double $\Upgamma^\Updelta$ of a generalized $n$-gon $\Upgamma$ of order $(s,s)$. By \cite[section 9]{part3}, it follows that if $\Upgamma^\Updelta$ and ${\Upgamma'}^\Updelta$ are isomorphic doubles of generalized $n$-gons $\Upgamma$ and $\Upgamma'$, 
then either $\Upgamma$ is isomorphic to $\Upgamma'$ or $\Upgamma$ is isomorphic to the point-line dual of $\Upgamma'$.

\ex
Let $n \geq 3$ be a positive integer. Let $\Upgamma$ be an ordinary $n$-gon. Then $\Upgamma^\Updelta$ is an ordinary $2n$-gon. 

\begin{figure}[h]
\begin{tikzpicture}[x=0.75pt,y=0.75pt,yscale=-1,xscale=1]

\draw  [line width=2.25]  (551,149) -- (512.5,215.68) -- (435.5,215.68) -- (397,149) -- (435.5,82.32) -- (512.5,82.32) -- cycle ;
\draw[black,fill=orange]  (551,149)  circle [x radius= 5, y radius= 5]   ;
\draw[black,fill=orange]  (512.5,215.68)  circle [x radius= 5, y radius= 5]   ;
\draw[black,fill=orange]  (435.5,215.68)  circle [x radius= 5, y radius= 5]   ;
\draw[black,fill=orange]  (397,149)  circle [x radius= 5, y radius= 5]   ;
\draw[black,fill=orange]  (435.5,82.32)  circle [x radius= 5, y radius= 5]   ;
\draw[black,fill=orange]  (512.5,82.32)  circle [x radius= 5, y radius= 5]   ;

\draw  [line width=2.25]  (196.76,193.34) -- (88.64,193.34) -- (142.7,99.7) -- cycle ;
\draw[black,fill=orange]  (196.76,193.34)  circle [x radius= 5, y radius= 5]   ;
\draw[black,fill=orange]   (88.64,193.34) circle [x radius= 5, y radius= 5]   ;
\draw[black,fill=orange]   (142.7,99.7)   circle [x radius= 5, y radius= 5]   ;

\draw   (280,143.5) node {$\scalebox{3}{$\rightarrow$}$} ;


\draw (144,78.4) node [anchor=north west][inner sep=0.75pt]    {$a$};
\draw (72,196.4) node [anchor=north west][inner sep=0.75pt]    {$b$};
\draw (198.76,196.74) node [anchor=north west][inner sep=0.75pt]    {$c$};
\draw (133,200.4) node [anchor=north west][inner sep=0.75pt]    {$A$};
\draw (178,129.4) node [anchor=north west][inner sep=0.75pt]    {$B$};
\draw (98,127.4) node [anchor=north west][inner sep=0.75pt]    {$C$};
\draw (514,59.4) node [anchor=north west][inner sep=0.75pt]    {$B$};
\draw (420,59.4) node [anchor=north west][inner sep=0.75pt]    {$a$};
\draw (561,142.4) node [anchor=north west][inner sep=0.75pt]    {$c$};
\draw (514.5,219.08) node [anchor=north west][inner sep=0.75pt]    {$A$};
\draw (421,219.4) node [anchor=north west][inner sep=0.75pt]    {$b$};
\draw (374,142.4) node [anchor=north west][inner sep=0.75pt]    {$C$};

\end{tikzpicture}
\caption{Doubling a triangle.}
\end{figure}
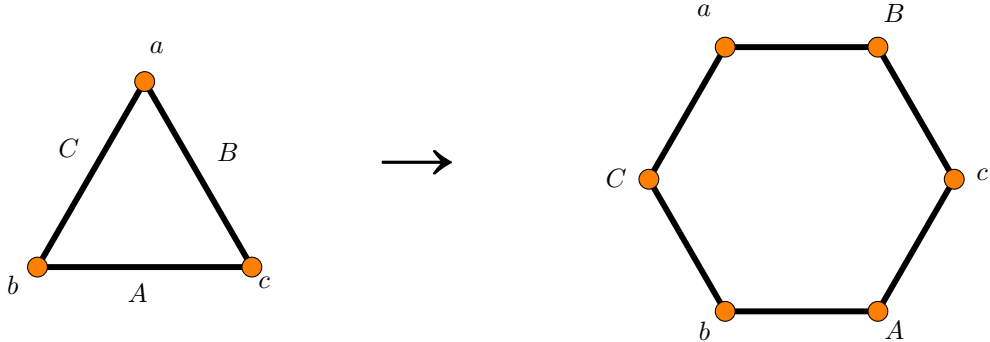

\subsection{Classical and nonclassical examples of thin generalized octagons}

Let $\F_q$ be a finite field. We work with homogeneous coordinates $(x_0 : x_1 : \cdots : x_m)$ in the projective 
space $\mathbb{P}^m(q)$. \\


Let $m = 4$. The $\F_q$-rational points and lines lying on a (parabolic) quadric with defining equation
\[ X_0X_1\ +\ X_2X_3\  + X_4^2\ =\ 0  \] 
define a thick generalized quadrangle $\mQ(4,q)$ with parameters $(q,q)$. \\

The {\em classical} generalized quadrangles of order $(s,s)$ are precisely the GQs $\mQ(4,q)$ with $q = s$ and the point-line duals of the latter. Doubling these 
quadrangles gives rise to classical thin generalized octagons. \\

Now let $\mO$ be a set of $q + 1$ points in $\mathbb{P}^2(q)$ such that no three of them are collinear --- by definition, this is an {\em oval} in $\mathbb{P}^2(q)$. Note that 
each point $x$ of $\mO$ is incident with a unique line which meets the oval only in $x$; this is a {\em tangent line}. 

We embed $\mathbb{P}^2(q)$ as a hyperplane in  $\mathbb{P}^3(q)$, and define a point-line geometry $\hT_2(\mO)$ as follows. 

\subsection*{Points} 

The points are defined as follows:  
\begin{itemize}
\item
A symbol $(\infty)$. 
\item
The points of $\mathbb{P}^3(q) \setminus \mathbb{P}^2(q)$. 
\item
Planes not contained in $\mathbb{P}^2(q)$ which meet $\mathbb{P}^2(q)$ in a tangent line to the oval.   
\end{itemize}

\subsection*{Lines} 

The lines come in two types:
\begin{itemize}
\item
The points of $\mO$. 
\item
The lines of $\mathbb{P}^3(q)$ which meet $\mathbb{P}^2(q)$ only in one point of $\mO$. 
\end{itemize}

It is easy to see that this geometry is a generalized quadrangle of order $(q,q)$, and it is known that $\hT_2(\mO) \cong \mQ(4,q)$ if and only if $\mO$ is a conic. If $q$ is even, there are many classes of examples of ovals which are not isomorphic to a conic, and these give rise to nonclassical quadrangles of order $(q,q)$, and hence to {\em nonclassical} generalized octagons of order $(q,1)$/$(1,q)$.

\subsection{Classical and nonclassical examples of thin generalized hexagons}

Besides the classical Desarguesian projective planes $\bP^2(q)$, where $q$ is any prime power, many infinite classes of finite projective planes are known which are not
isomorphic to a Desarguesian plane. Since the order of a plane is always of type $(N,N)$,  
such examples give rise to nonclassical generalized hexagons of order $(N,1)$/$(1,N)$.  

\subsection{Thin generalized \texorpdfstring{$12$}{12}-gons}

The only known finite generalized hexagons of order $(u,u)$ are, up to point-line duality, the {\em split Cayley hexagons} $\mathbf{H}(u)$ (in which case $u$ is a prime power). They give rise to classical thin generalized $12$-gons.

We refer to \cite[section 2]{POL} for more detailed information on classical polygons.

\subsection{Structural classification of finite \texorpdfstring{$\Fun$}{F1}-polygons} 

By recent work of Thas and Thas \cite{part3,part4} we can neatly describe all finite $\Fun$-generalized $m$-gons. The results below state first of all that all finite $\Fun$-polygons are proper, and secondly they describe the possible $\Fun$-structures which arise.  

\begin{theorem}[The planes \cite{part3}]
\label{GT}
Let $\Phi$ be an epimorphism of a thick projective plane $\mP$ onto a thin projective plane $\Updelta$ of order $(1,1)$.  Then exactly two classes of epimorphisms $\Phi$ occur (up to a suitable permutation of the points of $\Updelta$), and they are described as follows. 

\begin{itemize}
\item[{\rm(a)}] 
The points of $\Updelta$ are $\overline{a}, \overline{b}, \overline{c}$, with $\overline{a} \sim \overline{b} \sim \overline{c} \sim \overline{a}$, and put $\Phi^{-1} (\overline{x}) = \wt X$, with $\overline{x} \in \lbrace \overline{a}, \overline{b}, \overline{c} \rbrace$.

   Let $(\wt{A}, \wt{B})$, with $\wt{A} \ne \emptyset \ne \wt{B}$, be a partition of the set of all points incident with a line $L$ of $\mP$. Let  $\wt{C}$ consist of the points not incident with $L$. Furthermore, $\Phi^{-1}(\overline{a}\overline{b}) = \{ L\}$, $\Phi^{-1}(\overline{b}\overline{c})$ is the set of all lines distinct from $L$ but incident with a point of $\wt{B}$ and $\Phi^{-1}(\overline{a}\overline{c})$ is the set of all lines distinct from $L$ but incident with a point of $\wt{A}$.
   \item[{\rm(b)}] 
   The dual of $(a)$.
\end{itemize}
\end{theorem}

\begin{theorem}[The quadrangles \cite{part3}]
\label{JATGQ}
Let $\Phi$ be an epimorphism of a thick generalized quadrangle $\mathcal{S}$ of order $(s, t)$ onto a grid $\mathcal{G}$. Let $\mathcal{G}$ have order $(s^\prime, 1)$. Then $s^\prime = 1$ and exactly two classes of epimorphisms $\Phi$ occur (up to a suitable permutation of the points of $\mathcal{G})$.

\begin{itemize}
\item[{\rm(a)}] The points of $\mathcal{G}$ are $\overline{a}, \overline{b}, \overline{c}, \overline{d}$, with $\overline{a} \sim \overline{b} \sim \overline{c} \sim \overline{d} \sim \overline{a} $, and put $\Phi^{-1} (\overline{x}) = \wt X$, with $\overline{x} \in \lbrace \overline{a}, \overline{b}, \overline{c}, \overline{d} \rbrace$.

   Let $(\wt{A}, \wt{B})$, with $1 \le \vert \wt A \vert \le s, 1 \le \vert \wt B \vert \le s$, be a partition of the set of all points incident with a line $L$ of $\mathcal{S}$. Let  $\wt{C}$ consist of the points not incident with $L$ but collinear with a point of $\wt{B}$, and let $\wt{D}$ consist of the points not incident with $L$ but collinear with a point of $\wt{A}$. Further, $\Phi^{-1}(\overline{a}\overline{b}) = \{ L \}$, $\Phi^{-1}(\overline{b}\overline{c})$ is the set of all lines distinct from $L$ but incident with a point of $\wt{B}$, $\Phi^{-1}(\overline{a}\overline{d})$ is the set of all lines distinct from $L$ but incident with a point of $\wt{A}$ and $\Phi^{-1}(\overline{c}\overline{d})$ consists of all lines incident with at least one point of $\wt{C}$ and at least one point of $\wt{D}$.
\item[{\rm(b)}] The dual of $(a)$.
\end{itemize}
\end{theorem}

  \begin{theorem}[The hexagons \cite{part3}]
   \label{JATGH}
   Let $\Phi$ be an epimorphism of a thick generalized hexagon $\mS$ of order $(s, t)$ onto a thin generalized hexagon $\mG$ of order $(s^\prime, 1)$. Then $s^\prime = 1$ and exactly two classes of epimorphisms $\Phi$ occur (up to a suitable permutation of the points of $\mG$).
   
   \begin{itemize}
   \item[{\rm (a)}] The points of $\mG$ are $\ol{a}, \ol{b}, \ol{c}, \ol{d}, \ol{e}, \ol{f}$, with $\ol{a} \sim \ol{b} \sim \ol{c} \sim \ol{d} \sim \ol{e} \sim \ol{f} \sim \ol{a}$, and put $\Phi^{-1}(\ol{x}) = \wt{X}$, with $ \ol{x} \in \lbrace \ol{a}, \ol{b}, \ol{c}, \ol{d}, \ol{e}, \ol{f} \rbrace$.
   
   Let $(\wt{C}, \wt{B}), 1 \le \vert \wt{C} \vert \le s, 1 \le \vert \wt{B} \vert \le s$, be a partition of the set of all points incident with some line $L$ of $\mS$. Let $\wt{D}$ consist of the points not incident with $L$ but collinear with a point of $\wt{C}$, let $\wt{A}$ consist of the points not incident with $L$ but collinear with a point of $\wt{B}$, let $\wt{E}$ consist of the points not in $\wt{C} \cup \wt{D}$ but collinear with a point of $\wt{D}$, and let $\wt{F}$ consist of the points not in $\wt{A} \cup \wt{B}$ but collinear with a point of $\wt{A}$. Further, $\Phi^{-1}(\ol{b}\ol{c}) = \lbrace L \rbrace$, $\Phi^{-1}(\ol{c}\ol{d})$ is the set of all lines distinct from $L$ but incident with a point of $\wt{C}$, $\Phi^{-1}(\ol{a}\ol{b})$ is the set of all lines distinct from $L$ but incident with a point of $\wt{B}$, $\Phi^{-1}(\ol{d}\ol{e})$ is the set of all lines distinct from the lines of $\Phi^{-1}(\ol{d}\ol{c})$ but incident with a point of $\wt{D}$, $\Phi^{-1}(\ol{f}\ol{a})$ is the set of all lines distinct from the lines of $\Phi^{-1}(\ol{a}\ol{b})$ but incident with a point of $\wt{A}$, $\Phi^{-1}(\ol{d}\ol{e})$ is the set of all lines distinct from the lines 
   of $\Phi^{-1}(\ol{c}\ol{d})$ but incident with a point of $\widetilde{D}$,
   and $\Phi^{-1}(\ol{f}\ol{e})$ is the set of all lines not in $\Phi^{-1}(\ol{f}\ol{a})$ but incident with a point of $\wt{F}$ (that is, the set of all lines not in $\Phi^{-1}(\ol{e}\ol{d})$ but incident with a point of $\wt{E}).$
   
   \item[{\rm (b)}] The dual of $(a)$.
  
  \end{itemize}
  \end{theorem}
  
  \begin{theorem}[The octagons \cite{part4}]
\label{JATGO}
Let $\Phi$ be an epimorphism of a thick generalized octagon $\mS$ of order $(s, t)$ onto a thin generalized octagon $\mG$ of order $(s^\prime, 1)$. Then $s^\prime = 1$ and exactly two classes of epimorphisms $\Phi$ occur (up to a suitable permutation of the points of $\mG$).

\begin{itemize}
\item[{\rm (a)}] The points of $\mG$ are $\ol{a}, \ol{b}, \ol{c}, \ol{d}, \ol{e}, \ol{f}, \ol{g}, \ol{h}$, with $\ol{a} \sim \ol{b} \sim \ol{c} \sim \ol{d} \sim \ol{e} \sim \ol{f} \sim \ol{g} \sim \ol{h} \sim \ol{a}$, and put $\Phi^{-1}(\ol{x}) = \wt{X}$, with $\ol{x} \in \lbrace \ol{a}, \ol{b}, \ol{c}, \ol{d}, \ol{e}, \ol{f}, \ol{g}, \ol{h} \rbrace$.

Let $(\wt{C}, \wt{B}), 1 \le \vert \wt{C} \vert \le s, 1 \le \vert \wt{B} \vert \le s$, be a partition of the set of all points incident with a line $L$ of $\mS$. Let $\wt{D}$ consist of the points not incident with $L$ but collinear with a point of $\wt{C}$, let $\wt{A}$ consist of the points not incident with $L$ but collinear with a point of $\wt{B}$, let $\wt{E}$ consist of the points not in $\wt{C} \cup \wt{D}$ but collinear with a point of $\wt{D}$, let $\wt{H}$ consist of the points not in $\wt{A} \cup \wt{B}$ but collinear with a point of $\wt{A}$, let $\wt{F}$ consist of the points not in $\wt{D} \cup \wt{E}$ but collinear with a point of $\wt{E}$, and let $\wt{G}$ consist of the points not in $\wt{A} \cup \wt{H}$ but collinear with a point of $\wt{H}$.

Further, $\Phi^{-1}(\ol{c}\ol{b}) = \lbrace L \rbrace$, $\Phi^{-1}(\ol{c}\ol{d})$ is the set of all lines distinct from $L$ but incident with a point of $\wt{C}$, $\Phi^{-1}(\ol{a}\ol{b})$ is the set of all lines distinct from $L$ but incident with a point of $\wt{B}$, $\Phi^{-1}(\ol{d}\ol{e})$ is the set of all lines distinct from the lines of $\Phi^{-1}(\ol{d}\ol{c})$ but incident with a point of $\wt{D}$, $\Phi^{-1}(\ol{a}\ol{h})$ is the set of all lines distinct from the lines of $\Phi^{-1}(\ol{a}\ol{b})$ but incident with a point of $\wt{A}$, $\Phi^{-1}(\ol{e}\ol{f})$ is the set of all lines distinct from the lines of $\Phi^{-1}(\ol{d}\ol{e})$ but incident with a point of $\wt{E}$, $\Phi^{-1}(\ol{g}\ol{h})$ is the set of all lines distinct from the lines of $\Phi^{-1}(\ol{a}\ol{h})$ but incident with a point of $\wt{H}$, and $\Phi^{-1}(\ol{g}\ol{f})$ is the set of lines not in $\Phi^{-1}(\ol{h}\ol{g})$ but incident with a point of $\wt{G}$ (that is, the set of all lines not in $\Phi^{-1}(\ol{e}\ol{f})$ but incident with a point of $\wt{F}$).

\item[{\rm (b)}] The dual of $(a)$.
\end{itemize}
\end{theorem}

We summarize the results of this subsection as follows:
\begin{theorem}
Let $(\Upgamma,\Updelta,\upepsilon)$ be a finite $\Fun$-polygon with gonality $n \geq 3$. Then $(\Upgamma,\Updelta,\upepsilon)$ is a proper $\Fun$-polygon, and $\upepsilon: \Upgamma \mapsto \Updelta$ can be precisely described. 
\end{theorem} 

If $A = (\Upgamma,\Updelta,\upepsilon)$ and $B = (\Upgamma',\Updelta,\upepsilon)$ are two proper $\Fun$-generalized $n$-gons ($n \geq 3$), we can naturally define morphisms with source $A$ and target $B$ as morphisms of $n$-gons $\upgamma: \Upgamma \mapsto \Upgamma'$ for which the following diagram commutes:

\medskip

\begin{center}
{\scalebox{1.2}{
\begin{tikzpicture}[scale=2.2,text height=1.5ex, text depth=0.25ex]
\node (a0) at (0,3) {$\Updelta$};
\node (a1) [left=of a0] {$\Upgamma$};

\node (b1) [below=of a1] {$\Upgamma'$};

\draw[<<-,thick]
(a0) edge node[auto] {$\upepsilon'$} (b1)
(a0) edge node[auto] {$\upepsilon$} (a1);

\draw[->,orange,thick]
(a1) edge node[auto] {$\upgamma$} (b1);

\end{tikzpicture}
}}
\end{center}

 The morphisms $\upepsilon$ and $\upepsilon'$ endow $\Upgamma$ and $\Upgamma'$ with an {\em $\Fun$-structure.} The morphism $\upgamma$ is compatible with these structures.



\subsection{The infinite case}

When we allow the polygons to be infinite, there is no reasonable way to classify $\Fun$-structures, since free constructions are possible, cf. Thas and Thas \cite{part3, part4}.

\medskip
\section{Classification of \texorpdfstring{$\Fun$}{F1}-structures on \texorpdfstring{$\mathbb{P}^3(k)$}{P3(k)}, with \texorpdfstring{$k$}{k} a finite field}

In this section we classify $\Fun$-structures on $\mathbb{P}^3(k)$ with $k$ a finite field. We will make use of Theorem \ref{GT}, and will need the four types of $\Fun$-structures (in a finite plane $\mP$) described below, with ``base line $U$.'' In all the cases below, we consider an epimorphism 
\begin{equation}
\upepsilon: \ \mP \ \epiarrow{1.1cm} \ K(3),  
\end{equation}
where the vertices of $K(3)$ are called $a, b$ and $c$. Below, the point sets $A, B, C$ (which cover the points of $\mP$) will be respectively mapped to $a, b, c$. If the line $L$ of $\mP$ contains points of $A$ and $B$ (e. g.), then it will be mapped to the line $ab$ of $K(3)$. Each line of $\mP$ contains points from precisely two elements in $\{ A, B, C \}$. 

\begin{itemize}
\item[{\bf \fbox{Type I}}]
The point set of $U$ is partitioned in two distinct nonempty sets $A$ and $B$, and the points of $C$ are the points not on $U$. 
\item[\fbox{{\bf Type} $\widetilde{\mathbf{I}}$}]
The point set of $U$ is partitioned in two distinct nonempty sets $\wB$ and $\wC$, where $\vert \wB \vert = 1$. The single point of $\wB$ is incident with a line $V$ 
which is partitioned in sets $B$ and $A$ (where $\wB \subseteq B$).    
All the other points of  the plane are points of $C$. Note that this $\Fun$-structure is the same as that from the previous type, but where the line $U$ is different, relative to 
the structure.  
\item[\ {\bf \fbox{Type II}}]
We have that $A = \{ a \}$ is given by a single point, and 
the line set on $a$ is partitioned in two distinct nonempty sets $\widehat{B}$ and $\widehat{C}$; the points of $B$ consist of the points incident with the lines of $\widehat{B}$ except $a$,  and  the points of $C$ consist of the points incident with the lines of $\widehat{C}$ except again $a$. We also suppose that $\vert \widehat{B} \vert \geq 2$ and $\vert \widehat{C} \vert \geq 2$, to prevent overlap with the other types. 
\item[\ \fbox{{\bf Type} $\widetilde{\mathbf{II}}$}]
The line $U$ is partitioned in two point sets $\widehat{B}$ and $\widehat{C}$, where  $\vert \widehat{B} \vert \geq 2$ and $\vert \widehat{C} \vert \geq 2$ for the same reason as in the previous type. The set 
$A$ consists of one single point $a$ (not incident with $U$); now $x \in B$ if $x \in \widehat{B}$ or $xa \cap U$ is in $\widehat{B}$, and  $x \in C$ if $x \in \widehat{C}$ or $xa \cap U$ is in $\widehat{C}$. Note that this $\Fun$-structure is the same as that from the previous type, but where the line $U$ is different, relative to 
the structure.  

\end{itemize}

We now describe the possible $\Fun$-structures on $\mathbb{P}^3(k)$ relative to one fixed line. 

We first need the following proposition. 
\begin{proposition}[Dimension Theorem]
\label{propdim} 
Consider an epimorphism 
\begin{equation}
\upepsilon: \ \mP = \mathbb{P}^n(k) \ \epiarrow{1.1cm} \ K(n + 1),  
\end{equation}
where $n \geq 2$. Then we have that the image under $\upepsilon$ of a subspace of $\mP$ of dimension $m$, 
is again a subspace of dimension $m$ (that is, a sub-complete graph $K(m + 1)$). 
\end{proposition} 
\proof
First observe that the {\em image under $\upepsilon$ of a subspace is again a subspace (complete subgraph)}. For, let $\upbeta \ne \emptyset$ be a subspace of $\mathbb{P}^n(k)$. Let 
$x \ne y$ be points in $\upepsilon(\upbeta)$ (if those would not exists, $\upepsilon(\upbeta)$ consists of a point); Let $\upepsilon(\widehat{x}) = x$ and $\upepsilon(\widehat{y}) = y$, with 
$\widehat{x}$ and $\widehat{y}$ points in $\upbeta$. Then $\upepsilon(\widehat{x}\widehat{y}) = xy$ is a line in $\upepsilon(\upbeta)$. It follows that the latter is a subspace of $K(n + 1)$. \\

Next, let $\upbeta$ be a hyperplane of $\mathbb{P}^n(k)$. Suppose that $\upepsilon(\upbeta)$ is not a hyperplane of $K(n + 1)$, and let $x \in \mathbb{P}^n(k) \setminus \upbeta$ be a point so that $\upepsilon(x) \not\in \upepsilon(\upbeta)$ (and note that such points exist). Since $\upbeta$ is a hyperplane, every point of $\mathbb{P}^n(k)$ is on some line on $x$ and a point of $\upbeta$. It follows that $\upbeta(\mathbb{P}^n(k))$ is the complete graph on $\upepsilon(x)$ and the points of $\upepsilon(\upbeta)$, which is a contradiction since $\upepsilon$ is assumed to be surjective. So $\upepsilon(\upbeta)$ is a hyperplane in $K(n + 1)$.       

An easy induction argument now yields the theorem. 
\eop \\

We now proceed with the classification of $\Fun$-structures on $\mathbb{P}^3(k)$. so let 
\begin{equation}
\upepsilon: \ \mP = \mathbb{P}^3(k) \ \epiarrow{1.1cm} \ K(4)  
\end{equation}
be a surjection. 

Let $U$ be a fixed line in $\mathbb{P}^3(k)$, and consider the pencil of subplanes which contain $U$; if $\Pi$ is such a plane, then by Proposition \ref{propdim}, we know 
that $\upepsilon(\Pi)$ is a complete subgraph on $3$ points in $K(4)$. Relative to $U$, we know that $\Pi$ has to be of one of the types $\mathbf{I}$, $\widetilde{\mathbf{I}}$,   
$\mathbf{II}$, $\widetilde{\mathbf{II}}$.

\subsection*{The extremal case}

First we suppose that each line of $\mathbb{P}^3(k)$ meets one of $\{ A, B, C, D \}$ in precisely one point, so that the remaining points of the line are fully contained in another member 
of $\{ A, B, C, D \}$. 

Suppose that $U$ is a line for which $\vert U \cap A \vert = 1$, and let the remaining points be contained in $B$. Let $\Pi$ be a plane which contains $U$; then we have three 
possibilities for $\Pi$:
\begin{itemize}
\item[$\circled{1}$]
all points of $\Pi \setminus U$ are contained in $C$ or in $D$ (Type {\bf I});
\item[$\circled{2}$]
let $\{ a \}  = A \cap U$; there is one single line $V$ incident with $U \ne V$, all of whose points different from $a$ are in $C$ or $D$; the remaining points of $\Pi$ are contained in $B$ (Type $\widetilde{\mathbf{I}}$);
\item[$\circled{3}$]
there is one single line $V$ incident with $U \ne V$, and precisely one point incident with $V$ which is contained in $C$ or $D$; 
all other points on $V$ are in $A$; the remaining points of $\Pi$ are contained in $B$ (Type $\widetilde{\mathbf{I}}$).
\end{itemize}

By letting $\Pi$ vary, it is easy to see that only case $\circled{1}$ can occur (as otherwise lines exist which nontrivially meet three members of $\{ A, B, C, D \}$. Because of the extremal assumption, we now have a complete description of all possible extremal $\Fun$-structures which can be endowed on $\mathbb{P}^3(k)$, up to a permutation of $A, B, C, D$:
\begin{quote}
\begin{mdframed}[linewidth=2.5,linecolor=orange, topline=false,rightline=false,bottomline=false]
\item[$\circled{A}$]
{\em there is a line $U$ such that $\vert U \cap A \vert = 1$, the remaining points of $U$ are in $B$; furthermore, there is one plane $\Pi_C$ containing $U$ such that all point of 
$\Pi_C \setminus U$ are contained in $C$, and for all other planes $\Pi$ containing $U$, we have that all points of $\Pi \setminus U$ are contained in $D$.} 
\end{mdframed} 
\end{quote}

\subsection*{The non-extremal case}

Suppose now that we are not in the extremal case, so that there is a line $U$ which meets, example given, $A$ in more than one point and $B$ in more than one point. This means 
we can assume that any plane containing $U$ is either a plane of Type {\bf I} or {\bf II}.

First suppose that $U$ is contained in a plane $\Pi_C$ of Type {\bf II}, and suppose that the unique point in $\Pi_C$ which is not contained in $A \cup B$, is contained in $C$. Then it is easy to see that if there is at least one plane of Type {\bf I} containing $U$, then all planes on $U$ different from $\Pi_C$ must be of Type {\bf I} (since otherwise there are lines meeting three members  of $\{ A, B, C, D \}$ nontrivially); furthermore, clearly all points not on $U$ in such a plane are contained in $D$. 

Finally, suppose that all planes containing $U$ are of Type {\bf II}. Let $c \in C$ and $d \in D$, and note that these points are not incident with $U$; then $\widehat{U} := cd$ only contains points in $C \cup D$, and moreover, $C \cup D = \widehat{U}$. Suppose by way of contradiction that there is some line $V$ which meets at least three members of $\{ A, B, C, D \}$. Obviously, $V$ cannot meet $U$, and it must meet $\widehat{U}$. Consider the plane $\Pi := \Big\langle V, \widehat{U} \Big\rangle$. Then $\Pi$ meets $U$ in one point $u$ which we suppose to be in $A$ without loss of generality. It is now easy to see that all points of $\Pi \setminus V$ are contained in $A$, so that $V$ only meets two members of $\{ A, B, C, D \}$.

We end up with a complete description of all possible non-extremal $\Fun$-structures which can be endowed on $\mathbb{P}^3(k)$, up to a permutation of $A, B, C, D$:
\begin{quote}
\begin{mdframed}[linewidth=2.5,linecolor=orange, topline=false,rightline=false,bottomline=false]
{\em There is a line $U$ such that $\vert U \cap A \vert \geq 2$, and the remaining points of $U$ are in $B$ with $\vert B \vert \geq 1$, and we distinguish three cases.}
\item[$\circled{B}$]
{\em All planes containing $U$ are of Type {\bf I}; there is partition $\mC, \mD$ of the planes on $U$ such that if $\Pi \in \mC$, all points of $\Pi \setminus U$ are contained in $C$, and 
if $\Pi \in \mD$, all points of $\Pi \setminus U$ are in $D$;} 
\item[$\circled{C}$]
{\em There is one plane $\Pi_C$ of Type {\bf II} containing $U$, and its unique point not in $A \cup B$ is contained in $C$; for any other plane $\Pi$ on $U$ we have that the points 
of $\Pi \setminus U$ are contained in $D$;}
\item[$\circled{D}$] 
{\em There is a line $\widehat{U}$ for which $\vert C \cap \widehat{U} \vert \geq 1$, $\vert D \cap \widehat{U} \vert \geq 1$ and $C \cup D = \widehat{U}$; each point $u$ which is not 
incident with $U$ nor $\widehat{U}$ is incident with a unique line $W$ which meets $U$ and $\widehat{U}$, and if $V \cap U$ is a point of $A$, respectively $B$, then 
all points of $V \setminus V \cap \widehat{U}$ are points of $A$, respectively $B$.}\\
\end{mdframed} 
\end{quote}
It is easy to see that $\circled{B}$ is a special case of $\circled{D}$. And if we allow the case $\vert A \vert = 1$ in $\circled{B}$, then $\circled{A}$ is a special case of $\circled{B}$.
In short, we have the following classification.

\begin{theorem}
\label{Funstr}
Let $k$ be a finite field.
Consider the surjective morphism $\upepsilon: \ \mP = \mathbb{P}^3(k) \ \epiarrow{1cm} \ K(4)$; then the possible $\Fun$-structures endowed by $\upepsilon$ on $\mathbb{P}^3(k)$ are up to 
a permutation of $A, B, C, D$, described as follows.
 \begin{quote}
\begin{mdframed}[linewidth=2.5,linecolor=orange, topline=false,rightline=false,bottomline=false]
{\em There is a line $U$ such that $\vert U \cap A \vert \geq 1$ and $\vert U \cap B \vert \geq 1$, and $U \subseteq A \cup B$, and we distinguish two  cases.}
\item[{\rm $\circled{E}$}]
There is one plane $\Pi_C$ of Type {\bf II} containing $U$, and its unique point not in $A \cup B$ is contained in $C$; for any other plane $\Pi$ on $U$ we have that the points 
of $\Pi \setminus U$ are contained in $D$;
\item[{\rm $\circled{F}$}] 
There is a line $\widehat{U}$ for which $\vert C \cap \widehat{U} \vert \geq 1$, $\vert D \cap \widehat{U} \vert \geq 1$ and $C \cup D = \widehat{U}$; each point $u$ which is not 
incident with $U$ nor $\widehat{U}$ is incident with a unique line $W$ which meets $U$ and $\widehat{U}$, and if $V \cap U$ is a point of $A$, respectively $B$, then 
all points of $V \setminus V \cap \widehat{U}$ are points of $A$, respectively $B$.
\end{mdframed} 
\end{quote}
\eop
\end{theorem}

\begin{remark}{\rm 
Note that in contrast to Theorem \ref{GT}, cases $\circled{E}$ and $\circled{F}$ are self-dual. }
\end{remark}

\subsection*{The infinite case}

In the proof of Theorem \ref{Funstr} we have not used the finiteness condition at first sight. But in the proof of Theorem \ref{GT} of \cite{part3}, which we use at 
various points, the finiteness condition is crucial. Still, one might wonder whether Theorem \ref{GT} is also true for infinite planes as well. But it appears to be easy to 
construct counter examples. Here is one construction. Let $u, v, w, x$ be the four vertices of an ordinary $4$-gon $\Upgamma$ (given in cyclic order), and add two different points per side. Also, add the line 
$xv$ (with no extra points).  Let $\Updelta$ be an ordinary triangle with vertices $a, b, c$.  
Now define $\upepsilon$ as follows:
\begin{equation}
\begin{cases}
\upepsilon(u) &= \upepsilon(w) = a \\
\upepsilon(v) &= c \\
\upepsilon(x) &= b
\end{cases}
\end{equation}

For each line $L$ of $\Upgamma$, map the remaining two points of $L$ which do not have an image yet, to two different vertices of $\Updelta$ (taking into account that 
$L$ meets precisely two sets of $\upepsilon^{-1}(a)$,  $\upepsilon^{-1}(b)$, $\upepsilon^{-1}(c)$). Then $\upepsilon$ is a surjective morphism of incidence geometries, and by \cite[section 6]{part4}, 
$\upepsilon$ can be extended to a surjective morphism 
\begin{equation}
\overline{\upepsilon}:\ \overline{\Upgamma}\ \epiarrow{1.1cm}\ \Updelta,
\end{equation}
where $\overline{\Upgamma}$ is the free projective plane generated by $\Upgamma$. Clearly, the $\Fun$-structure which arises is not described by Theorem \ref{GT}. \\


\end{document}